

\input epsf.tex

\def\2{{1\over 2}}

\def\d{\delta}
\def\a{\alpha}
\def\b{\beta}
\def\g{\gamma}

\def\e{\epsilon}
\def\l{\lambda}

\def\D{\Delta}

\def\fun#1#2#3{#1\colon #2\rightarrow #3}

\def\frac#1#2{{{#1} \over {#2}}}

\def\st{\;\colon\;}
\def\tends{\rightarrow}

\def\dr{ {\rm d} }

\def\R{{\bf R}}
\def\N{{\bf N}}
\def\Z{{\bf Z}}

\def\T{{\bf T}}

\def\thm#1{\vskip 1 pc\noindent{\bf Theorem #1.\quad}\sl}
\def\lem#1{\vskip 1 pc\noindent{\bf Lemma #1.\quad}\sl}
\def\prop#1{\vskip 1 pc\noindent{\bf Proposition #1.\quad}\sl}
\def\cor#1{\vskip 1 pc\noindent{\bf Corollary #1.\quad}\sl}

\def\proof{\rm\vskip 1 pc\noindent{\bf Proof.\quad}}
\def\fin{\par\hfill $\backslash\backslash\backslash$\vskip 1 pc}
\def\txt#1{\quad\hbox{#1}\quad}
\def\m{\mu}
\def\L{{\cal L}}

\def\G{\Gamma}

\def\r{\rho}

\def\2{\frac{1}{2}}
\def\inn#1#2{{\langle #1 ,#2\rangle}}

\def\part{{\partial_{x}}}
\def\div{{\rm div}}

\def\pprime{{{}^\prime{}^\prime}}
\def\ec{{\cal E}}

\def\fc{{\cal F}}

\def\dc{{\cal D}}
\def\pc{{\cal P}}
\def\vc{{\cal V}}

\def\zc{{\cal Z}}
\def\tc{{\cal T}}
\def\ent{{\rm Ent}}



\baselineskip= 17.2pt plus 0.6pt
\font\titlefont=cmr17
\centerline{\titlefont An entropy generation formula}
\vskip 1 pc
\centerline{\titlefont on $RCD(K,\infty)$ spaces}
\vskip 4pc
\font\titlefont=cmr12
\centerline{         \titlefont {Ugo Bessi}\footnote*{{\rm 
Dipartimento di Matematica, Universit\`a\ Roma Tre, Largo S. 
Leonardo Murialdo, 00146 Roma, Italy.}}   }{}\footnote{}{
{{\tt email:} {\tt bessi@matrm3.mat.uniroma3.it} Work partially supported by the PRIN2009 grant "Critical Point Theory and Perturbative Methods for Nonlinear Differential Equations}} 
\vskip 0.5 pc
 
\par
\vskip 2pc
\centerline{\bf Abstract}

J. Feng and T. Nguyen have shown that the solutions of the Fokker-Planck equation in 
$\R^d$ satisfy an entropy generation formula. We prove that, in compact metric measure spaces with the $RCD(K,\infty)$ property, a similar result holds for curves of measures whose density is bounded away from zero and infinity. We use this fact to show the existence of minimal characteristics for the stochastic value function.

\vskip 2 pc
\centerline{\bf  Introduction}
\vskip 1 pc

Let $\T^d=\frac{\R^d}{\Z^d}$ denote the $d$-dimensional torus and let $\pc(\T^d)$ denote the set of the Borel probability measures on $\T^d$. Let the two functions  
$\fun{F}{(-\infty,0]\times\T^d}{\R}$ and $\fun{U}{\pc(\T^d)}{\R}$ be continuous. The stochastic value function 
$$\fun{U}{(-\infty,0]\times\pc(\T^d)}{\R}$$
is defined as 
$$U(t,\eta)=\inf\left\{
\int_t^0\dr s\int_{\T^d}\left[
\2|X(s,x)|^2-F(s,x)
\right]    \dr\mu_s(x)+   U(\mu_0)
\right\}   \eqno (1)$$
where $t\le 0$ and $\fun{\mu}{[t,0]}{\pc(\T^d)}$ is a weak solution of the Fokker-Planck equation
$$\left\{
\matrix{
&\partial_s\mu_s-\frac{\b}{2}\D\mu_s+\div(X_s\mu_s)=0,\quad s\in[t,0]\cr
&\mu_t=\eta    .    
}
\right.   \eqno (2)$$
The $\inf$ in (1) is over all "reasonable" vector fields $X$. 

A natural question is whether the $\inf$ in (1) is a minimum. When $U$ is linear, say
$$U(\mu)=\int_{\T^d}f(x)\dr\mu(x)$$
for some $f\in C(\T^d)$, existence of minimisers was proven by Fleming in [12], essentially by establishing the optimal drift $X$; it is an approach that carries over, word for word, to 
$RCD(K,\infty)$ spaces ([8]). 

When $U$ is not linear, existence of minimisers was proven much more recently by J. Feng and T. Nguyen in [11]; the following two facts are at the core of their proof. 

\noindent 1) Let $\fc$ be a family of solutions $(\mu,X)$ of (2) which satisfies 
$$\sup_{(\mu,X)\in\fc}
\int_t^0\dr s\int_{\T^d}|X(s,x)|^2\dr\mu_s(x)<+\infty    $$
and let
$$E=\{
\mu\st(\mu,X)\in\fc\txt{for some}X
\}  .  $$
Then, $E$ is relatively compact in $C([t,0],\pc(\T^d))$. 

\vskip 1pc

\noindent 2) The drift in (2) is not determined uniquely by the curve $\mu$; however, if we associate to $\mu$ the drift $V_\mu$ of minimal norm in $L^2([t,0]\times\T^d,\L^1\otimes\mu_s)$, then the map from $C([t,0],\pc(\T^d))$ to $\R$ given by 
$$\fun{}{\mu}{
\int_t^0\dr s\int_{\T^d}|V_\mu(s,x)|^2\dr\mu_s(x)
}  $$
is lower semicontinuous.

The aim of this paper is to see whether this result of [11] can be generalised to a particular class of metric measure spaces, those with the $RCD(K,\infty)$ property. A metric measure space $(S,d,m)$ is simply a metric space $(S,d)$ with a Borel measure $m$ attached. We shall suppose that $(S,d)$ is compact and that $m$ is a probability measure; we point out at the outset that this is already a major simplification with respect to [11], which is set in the unbounded $\R^d$ for the Lebesgue measure. As noted at page 333 of [11], the bounded case can be treated by a standard mollification technique, and this is exactly what we do. We shall mollify $\mu_t$ in the usual way, i. e. applying the heat flow to it; the main reason we are working on $RCD(K,\infty)$ spaces is that their heat flow is very well behaved. 

We briefly explain how the Fokker-Plank equation translates to $RCD(K,\infty)$ spaces; we follow [5] and [13], which generalise the continuity equation to a very general class of metric measure spaces; in section 2 below, we shall briefly retrace the history of the problem. 

Starting from the torus, we recall that 
$\fun{\mu}{[t,0]}{\pc(\T^d)}$ is a weak solution of (2) if 
$$\int_t^0\dr s\int_{\T^d}
\left[
\partial_s\phi(s,x)+\frac{\b}{2}\D\phi(s,x)+\inn{\nabla\phi(s,x)}{X(s,x)}
\right]   \dr\mu_s(x)=0\qquad
\forall\phi\in C^\infty_0((t,0)\times\T^d)  .  \eqno (3)$$

First of all, we need a Laplacian; we recall from [6] that, if $(S,d,m)$ has the 
$RCD(K,\infty)$ property, then it admits an operator $\D_\ec$ with many of the properties of the standard Laplacian on $\T^d$: linearity, the integration by parts formula, a linear heat flow... 
Also the Dirichlet integral 
$$\2\int_{\T^d}|\nabla u|^2\dr x$$
has a counterpart; it is called Cheeger's energy, and usually denoted with $\2\ec$. Moreover, there is a "carr\'e\ de champs" operator $\G(\phi,\psi)$ which shares many properties with $\inn{\nabla\phi}{\nabla\psi}$.

Second, we need a class $\tc$ of test functions; these are readily available because, as proven in [7], the functions with bounded carr\'e\ de champs and bounded Laplacian are dense in the domain of Cheeger's energy. 

Next, we need a drift. Following [11] and [13], the right way to look at the drift is as an operator on test functions, bounded in a suitable norm. Following the approach of [5], we  consider the seminorm on the test functions $\tc$  
$$||\phi||^2_{\vc(\mu)}\colon=
\int_t^0\dr s\int_S\G(\phi_s,\phi_s)\dr\mu_s(x)  .  $$
After identifying $u$ and $v$ if $||u-v||_{\vc(\mu)}=0$, we can define $\vc(\mu)$ as the completion of $\tc$ with respect to $||\cdot||_{\vc(\mu)}$; in section 2 below we shall see that $\vc(\mu)$ is a Hilbert space. We look at the drift as the operator 
$$\fun{\tilde X}{\phi}{
-\int_t^0\dr s\int_S\left[
\partial_s\phi_s+\frac{\b}{2}\D_\ec\phi_s
\right]  \dr\mu_s
}   . $$
If $\tilde X$ is bounded on $\vc(\mu)$, i. e. if 
$$||
\dot\mu-\frac{\b}{2}\D_\ec\mu
||_{\vc^\prime(\mu)}\colon=$$
$$\sup\left\{
-\int_t^0\dr s\int_S\left[
\partial_s\phi_s+\frac{\b}{2}\D_\ec\phi_s
\right]   \dr\mu_s   \st
\phi\in\tc\txt{and}
\int_t^0\dr s\int_S\G(\phi_s,\phi_s)\dr\mu_s(x)\le 1
\right\}   <+\infty  \eqno (4)$$
then we can apply the Riesz representation theorem and find $X\in\vc(\mu)$ such that 
$$\int_t^0\dr s\int_S\left[
\partial_s\phi_s+\frac{\b}{2}\D_\ec\phi_s
\right]  \dr\mu_s+\inn{X}{\phi}_{\vc(\mu)}=0\qquad
\forall\phi\in\tc  .  $$
This will be our weak form of the Fokker-Planck equation. 

This way of defining the drift as a linear operator has the further advantage that the drift is automatically unique. 

The last thing we need, and only for the proof of corollary 1 below, is a very strong property of $(S,d,m)$, the so-called $L^1 -L^\infty$ regularisation property. 

We want to prove the following theorem; in section 1 we shall give definitions and references for all the terms which appear in it. 

\thm{1} Let $(S,d,m)$ be a $RCD(K,\infty)$ space; we suppose that $(S,d)$ is compact and that $m$ is a probability measure. Let $t<0$ and let $\fun{\mu}{[t,0]}{\pc(S)}$ be a continuous curve such that $\mu_s=\r_s m$ for all $s\in[t,0]$. Let us suppose that 
$$||
\dot\mu-\frac{\b}{2}\D_\ec\mu
||_{\vc^\prime(\mu)}<+\infty  .  $$
We also suppose that there is $C>0$ such that, for all $s\in[t,0]$,
$$\frac{1}{C}\le \r_s(x)\le C\txt{for $m$ a. e.} x\in S  .  \eqno (5)$$
Then,
$$||
\dot\mu-\frac{\b}{2}\D_\ec\mu
||_{\vc^\prime(\mu)}^2\ge$$
$$\int_t^0||\dot\mu_s||^2\dr s+
\b\ent_m(\mu_0)-\b\ent_m(\mu_t)  , \eqno (6)$$
where we have denoted by $||\dot\mu_s||$ the metric derivative of $\mu$ in $\pc(S)$ with respect to the 2-Wasserstein distance. 

\rm

\vskip 1pc

The following corollary yields the existence of minimisers for the value function.  

\cor{1} Let $(S,d,m)$ be a compact $RCD(K,\infty)$; we suppose that $(S,d)$ is compact and that $m$ is a probability measure. Then, the following holds. 

\noindent 1) Let $\mu^n\tends \mu$ in $C([t,0],\pc(S))$; let $\mu^n_s=\r^n_sm$ and let us suppose that 
$$\{ \r^n_s\st n\ge 1,\quad s\in[t,0] \}  $$
is uniformly integrable. Then, 
$$||\dot\mu-\frac{\b}{2}\D_\e\mu
||_{\vc^\prime(\mu)}\le
\liminf_{n\tends+\infty}
||\dot\mu^n-\frac{\b}{2}\D_\e\mu^n
||_{\vc^\prime(\mu^n)}  .  $$ 

\noindent 2) If $\fc\subset C([t,0],\pc(S))$ satisfies
$$M\colon=\sup_{\mu\in\fc}||
\dot\mu-\frac{\b}{2}\D_\ec\m
||_{\vc^\prime(\mu)}   <+\infty   ,  \eqno (7)$$
$$\sup_{\mu\in\fc}\ent_m(\mu_t)<+\infty   \eqno (8)$$
and if there is $c(\mu)>0$ such that $\mu_t=\r_tm$ with
$$\frac{1}{c(\mu)}\le\r_s(x)\le  c(\mu)\txt{a. e.}x\in S,\quad\forall s\in[t,0]
\eqno (9)$$
then $\fc$ is relatively compact in $C([t,0],\pc(S))$.

\noindent 3) Let $F\in C((-\infty,0]\times S,\R)$, let $U\in C(\pc(S),\R)$, let $t<0$ and 
$\nu\in\pc(S)$ be such that $\ent_m(\nu)<+\infty$. Let us suppose moreover that the heat kernel of $(S,d,m)$ satisfies formula (1.8) below. Then, the following infimum is attained.
$$U(t,\nu)\colon=\inf\Big\{
\2||\dot\mu-\frac{\b}{2}\D_\ec\mu||_{\vc^\prime(\mu)}^2+
\int_t^0\dr s\int_S F(s,x)\dr\mu_s(x)+
U(\mu_0)\st  $$
$$\mu\in C([t,0],\pc(S))\txt{and}\mu_t=\nu
\Big\}   .    \eqno (10)  $$

\rm

\vskip 1pc

The paper is organised as follows: in section 1 we recall some of the properties of 
$RCD(K,\infty)$ spaces from [2], [4], [6] and [7]; in section 2, we introduce the class $\tc$ of test functions; after an exposition of some results of [5], we shall arrive to the notion of weak solutions of the continuity and Fokker-Planck equations. In section 3, we review the standard method to regularise the solutions of Fokker-Planck. In section 4, we prove (6) when the curve $\mu$ is sufficiently regular. In section 5, we prove (6) in the general case, approximating the drift $X$ with the regularised drifts of section 3. In section 6, we prove corollary 1.  

\vskip 1pc

\noindent{\bf Acknowledgements.} The author would like to thank the referees for the patience and the useful suggestions.

\vskip 2 pc
\centerline{\bf  \S 1}
\centerline{\bf Preliminaries and notation}
\vskip 1 pc

\noindent{\bf The metric space.} We shall suppose that $(S,d)$ is a compact metric space with the geodesic property: in other words, for all $x,y\in S$ there is $\fun{\g}{[0,1]}{S}$ such that $\g_0=x$, $\g_1=y$ and, for all $0\le s\le t\le 1$, 
$$d(\g_s,\g_t)=(t-s)d(x,y)  .  $$
Such a curve is called a constant speed geodesic connecting $x$ and $y$. 

If $\fun{f}{S}{\R}$ is a function, we define its slope as 
$$|Df|(x)=\limsup_{y\tends x}\frac{|f(y)-f(x)|}{d(x,y)}  .  \eqno (1.1)$$

\noindent{\bf Absolutely continuous curves.} We say that $\fun{\g}{[a,b]}{S}$ belongs to 
$AC^2([a,b],S)$ if there is $\a\in L^2(a,b)$ such that, for all $a\le s\le t\le b$ we have that
$$d(\g_s,\g_t)\le\int_s^t\a(\tau)\dr\tau  .  $$
By [3], there is a function $\a$ of minimal $L^2$ norm for which the formula above holds; we call it $||\dot\g_\tau||$. 

An obvious remark is that constant speed geodesics are an example of absolutely continuous curves. 

\noindent{\bf The Wasserstein distance.} Let $\pc(S)$ denote the set of all Borel probability measures on $S$; for $\mu,\nu\in\pc(S)$ we set 
$$W_2^2(\mu,\nu)=\min
\int_{S\times S}d^2(x,y)\dr\g(x,y)$$
where the minimum is over all Borel probability measures $\g$ on $S\times S$ whose first and second marginals are $\mu$ and $\nu$ respectively. 

It can be proven ([17], [3], [1]) that the minimum is attained and that $W_2$ induces the weak$\ast$ topology on $\pc(S)$. 

\noindent{\bf Cheeger's energy and its flow.} From now on we shall fix one $m\in\pc(S)$ and consider the triple $(S,d,m)$. As usual, we denote by $F_\sharp\mu$ the push-forward of a measure $\mu$ by a map $F$; for $s\in[a,b]$ we denote by 
$$\fun{e_s}{C([a,b],S)}{S},\qquad
\fun{e_s}{\g}{\g_s} $$
the evaluation map. 

Let $\pi$ be a Borel probability measure on $C([0,1],S)$; we say that $\pi$ is a test plan if the following points hold. 

\noindent 1) $\pi$ concentrates on $AC^2([0,1],S)$.

\noindent 2) We have that
$$\int_{C([0,1],S)}\dr\pi(\g)\int_0^1||\dot\g_\tau||^2\dr\tau<+\infty  .  $$

\noindent 3) Let us define $\mu_s=(e_s)_\sharp\pi$; we ask that  $\mu_s=\r_s m$ for all 
$s\in[0,1]$ and that there is $C>0$ such that 
$$||\r_s||_{L^\infty(S,m)}\le C
\txt{for all} s\in[0,1]  .  $$

Let $\fun{f}{S}{\R}$ be a function; we say that a Borel function $\fun{g}{S}{[0,+\infty]}$ is a weak upper slope of $f$ if, given any test plan $\pi$, 
$$|f(\g_1)-f(\g_0)|\le\int_0^1g(\g_\tau)||\dot\g_\tau||\dr\tau$$
holds for for $\pi$ a. e. curve. 

If $f\in L^2(S,m)$ has a weak upper slope $g\in L^2(S,m)$, then it has an upper slope 
$|Df|_w$ minimal in the following two senses: its $L^2$ norm is smaller than all the other upper slopes and 
$$|Df|_w(x)\le h(x)  \txt{for $m$ a. e. $x\in S$}$$
if $h$ is an upper slope. 

Naturally, $|Df|_w$ is not unrelated to $|Df|$ defined in (1.1); this relationship is, roughly, Cheeger's definition of $|Df|_w$, for which we refer the reader to [10] and [4]. 

The double of Cheeger's energy 
$$\fun{\ec}{L^2(S,m)}{[0,+\infty]}$$
is defined as 
$$\ec(f)=\int_S|Df|_w^2\dr m   \eqno (1.2)$$
if $f$ has a weak upper slope in $L^2(S,m)$, and $\ec(f)=+\infty$ otherwise. 

Though on the right hand side of (1.2) there is a square, $\ec(f)$ in general is not a quadratic form (see remark 4.7 of [4] for an example). However, it is a convex functional; in particular, we can fix once and for all a positive constant $\b$ and define a gradient flow
$$\dot f_\tau\in-\frac{\b}{4}\partial\ec(f_t)$$
starting from every $f_0\in L^2(S,m)$. Since the standard heat flow is the gradient flow of the Dirichlet energy, we can think that $f_t$ is a kind of heat flow in metric measure spaces. 

We also have a kind of Laplacian: if $\partial\ec(u)\not=\emptyset$, then we define 
$-\D_\ec u$ as the element of smallest norm in $\partial\ec(u)$. Since $\ec$ is not necessarily quadratic, $\D_\ec$ and the heat flow are not necessarily linear. 

\noindent{\bf The entropy functional.} We define
$$\fun{\ent_m}{\pc(S)}{[0,+\infty]}$$
$$\ent_m(\mu)=\left\{
\eqalign{
\int_S\r\log\r\dr m &\txt{if}\mu=\r m\cr
+\infty  &\txt{otherwise.}
}
\right.   $$
It turns out ([4]) that, if $\ent_m$ satisfies the $CD(K,\infty)$ condition, which we state below, then each $\mu\in\pc(S)$ with $\ent_m(\mu)<+\infty$ is the starting point of a unique gradient flow (in the EDE sense, for whose definition we refer to [4]) of 
$\frac{\b}{2}\ent_m$. 

Moreover, if $\mu=\r_0m$ and $\r_0$ is in $L^2(S,m)$, then the gradient flow of 
$\frac{\b}{2}\ent_m$ has the form $\mu_\tau=\r_\tau m$, where $\r_\tau$ is the gradient flow of $\frac{\b}{4}\ec$ starting at $\r_0$. In the following, we shall set 
$H_\tau\mu=\mu_\tau$, where $\mu_\tau$ is the gradient flow of $\frac{\b}{2}\ent_m$ starting at $\mu$; this defines a semigroup $H_\tau$ on the domain of $\ent_m$. 

\noindent{\bf The $RCD(K,\infty)$ condition.} Following [4] and [2], we say that $(S,d,m)$ is a $RCD(K,\infty)$ space if the following two conditions hold.

\noindent 1) $\ec$ is quadratic, i. e. setting
$$\dc(\ec)=\{
u\in L^2(S,m)\st \ec(u)<+\infty
\}   $$
the parallelogram equality holds 
$$\ec(u)+\ec(v)=\2[
\ec(u+v)+\ec(u-v)
]  \qquad\forall u,v\in\dc(\ec)   .  \eqno (1.3)$$

\noindent 2) The $CD(K,\infty)$ condition holds; in other words, for all 
$\tilde\mu_0,\tilde\mu_1\in\pc(S)$ with $\ent_m(\tilde\mu_0),\ent_m(\tilde\mu_1)<+\infty$, there is a constant speed geodesic $\mu_t$ such that $\mu_0=\tilde\mu_0$, 
$\mu_1=\tilde\mu_1$ and 
$$\ent_m(\mu_t)\le
(1-t)\ent_m\mu_0+t\ent_m\mu_1-
\frac{K}{2}t(1-t)W_2^2(\mu_0,\mu_1) \qquad\forall t\in[0,1]  .  $$

\vskip 1pc

We recall one consequence of the $RCD(K,\infty)$ condition: if $\mu,\nu$ are in the domain of $\ent_m$, then 
$$W_2(H_s\mu,H_s\nu)\le e^{-\frac{K\b}{2}s} W_2(\mu,\nu)   \eqno (1.4)$$
for the same $K$ in the definition of $RCD(K,\infty)$. 

As we saw above, there is a gradient flow of $\ent_m$ starting from every $\mu\in\pc(S)$ with $\ent_m(\mu)<+\infty$; since it is easy to see that measures with bounded entropy are dense in $\pc(S)$,  the uniform continuity provided by (1.4) allows us to extend 
$H_s$ to a semigroup defined on the whole of $\pc(S)$. 

\vskip 1pc

\noindent{\bf The carr\'e\ de champs.} From now on we shall suppose that $(S,d,m)$ is a $RCD(K,\infty)$ space; as a consequence of (1.3) we can define a symmetric bilinear form 
$$\fun{\ec}{\dc(\ec)\times\dc(\ec)}{\R}$$
by
$$\ec(u,v)=\frac{1}{4}[
\ec(u+v)-\ec(u-v)
]  .  $$
The important fact is that $|Du|_w^2$ is a bilinear form too. Namely ([6], [7]) there is a symmetric bilinear operator
$$\fun{\G}{\dc(\ec)\times\dc(\ec)}{L^1(S,m)}  $$
such that
$$\ec(u,v)=\int_S\G(u,v)\dr m\qquad
\forall(u,v)\in\dc(\ec)\times\dc(\ec)  .  $$
The operator $\G$ is called the carr\'e\ de champs of $\ec$. We recall that  
$\dc(\D_\ec)\subset\dc(\ec)$ and that the integration by parts formula holds
$$\inn{-\D_\ec u}{v}_{L^2(S,m)}=\ec(u,v)=
\int_S\G(u,v)\dr m\qquad
\forall u\in \dc(\D_\ec), \qquad \forall v\in\dc(\ec) .   \eqno (1.5)$$

\noindent{\bf The heat semigroup.} The operator $\frac{-\b}{2}\D_\ec$ is monotone maximal; thus (see for instance [9]) it gives rise to a semigroup, backward in time,
$$\fun{P_{-s}}{L^2(S,m)}{L^2(S,m)}\qquad
P_{-s}=e^{
s\frac{\b}{2}\D_\ec
}   \txt{for $s\ge 0$.}  $$
Now $P_{-s}$ is induced by a Brownian motion; namely, there is a probability measure 
${\bf P}^{-s,x}$ on $C([-s,+\infty],S)$ which concentrates on the curves $\g$ with 
$\g(-s)=x$; denoting by $E^{-s,x}$ the expectation with respect to 
${\bf P}^{-s,x}$, we have that
$$(P_{-s}f)(x)=E^{-s,x}[f(\g_0)]  .  $$
By the way, the existence of the Brownian motion is one of the reasons for point 1) in the definition of $RCD(K,\infty)$ spaces: we want $\frac{\b}{2}\ec$ to be the Dirichlet form of the Brownian motion, and Dirichlet forms are quadratic. 

Since $m(S)=1$ and $S$ is compact, we have that $C(S)\subset L^2(S,m)$; essentially by (1.4), we have that $\fun{P_{-s}}{C(S)}{C(S)}$; thus, $P_{-s}$ has an adjoint on the space of Borel measures on $S$. This adjoint is the operator $H_s$ we defined above: 
$$\int_S (P_{-s}f)\dr\mu=
\int_S f\dr(H_s\mu)  \qquad\forall f\in C(S) 
\quad\forall\mu\in\pc(S) .  $$

As we noted above, if $\mu=\r m$ with $\r\in L^2(S,m)$, we have that 
$$H_s(\mu)=(P_{-s}\r)m  .  \eqno (1.6)$$

We recall that $P_{-s}$ has a symmetric kernel; namely,  for all $s>0$ there is a measurable function 
$\fun{p_{-s}}{S\times S}{[0,+\infty)}$ such that

\noindent 1) $p_{-s}(x,y)=p_{-s}(y,x)$ for all $s>0$ and $m\otimes m$ a. e. 
$(x,y)\in S\times S$, and 

\noindent 2) for all $f\in L^2(S,m)$ and all $s>0$, for $m$ a. e. $x\in S$, we have
$$(P_{-s}f)(x)=\int_S f(y)p_{-s}(x,y)\dr m(y)  .  \eqno (1.7)$$

For corollary 1 we shall need a very strong property of $p_{-s}$, namely that for all $s>0$ there is $a_s>0$ such that 
$$\frac{1}{a_s}\le p_{-s}(x,y)\le a_s
\txt{for $m\otimes m$ a. e. $(x,y)\in S\times S$.}   \eqno (1.8)$$
This implies the $L^1-L^\infty$ regularisation property; indeed, by (1.7) and (1.8), 
$$\frac{1}{a_s}||f||_{L^1(S,m)}\le
 P_{-s}(f)\le a_s||f||_{L^1(S,m)}  .  \eqno (1.9)$$

\noindent{\bf The smooth functions.} We recall (see for instance section 3 of [15]) that, since $\ec$ is lower semicontinuous for the $L^2$ topology, $\dc(\ec)$ with the inner  product defined by
$$\inn{u}{v}_{\dc(\ec)}=\inn{u}{v}_{L^2(S,m)}+\ec(u,v)$$
is a Hilbert space; we define $||u||_{\dc(\ec)}^2=\inn{u}{u}_{\dc(\ec)}$ as the natural norm.

We recall from [7] that the subspaces 
$$V^1_\infty=\{
u\in\dc(\ec)\cap L^\infty(S,m)\st \G(u,u)\in L^\infty(S,m)
\}   $$
and 
$$V^2_\infty=\{
u\in V^1_\infty\st\D_\ec u\in L^\infty(S,m)
\}  $$
are dense in $\dc(\ec)$. We note that $V^1_\infty$ with the norm
$$||u||_{V^1_\infty}=||u||_{L^\infty(S,m)}+||\G(u,u)^\2||_{L^\infty(S,m)}$$
is a Banach space. We briefly prove that it is complete. Let $\{ u_n \}_{n\ge 1}$ be a Cauchy sequence in $V^1_\infty$: given $\e>0$ we can find $N\in\N$ such that
$$||u_n-u_l||_\infty+||\G^\2(u_n-u_l,u_n-u_l)||_\infty\le\e
\txt{if}n,l\ge N  .  \eqno (1.10)$$
By the formula above we get that $\{ u_n \}_{n\ge 1}$ is Cauchy also in $\dc(\ec)$; since 
$\dc(\ec)$ is complete, we have that $u_n\tends\bar u$ in $\dc(\ec)$. 

Recall that $\G$ is continuous from $\dc(\ec)\times\dc(\ec)$ into $L^1$ and 
$u_l\tends\bar u$ in $\dc(\ec)$; thus, we can fix $n\in\N$, let $l\tends+\infty$ and get that, for all $n$, 
$$u_n-u_l\tends u_n-\bar u\txt{in $L^2$ and}$$
$$\G(u_n-u_l,u_n-u_l)\tends\G(u_n-\bar u,u_n-\bar u)
\txt{in}L^1  .  $$
Up to subsequences, the convergence in the two formulas above is $m$-a.e.; together with (1.10), this implies that 
$$||u_n-\bar u||_\infty+||\G^\2(u_n-\bar u,u_n-\bar u)||_\infty\le\e
\txt{if}n\ge N    $$
i. e. that $\bar u\in V^1_\infty$ and  $u_n\tends\bar u$ in $V^1_\infty$. 

With the same argument, it is possible to show that that $V^2_\infty$ with the norm 
$$||u||_{V^2_\infty}\colon=||u||_{V^1_\infty}+||\D_\ec u||_{L^\infty}$$
is a Banach space. 

We shall need the following formula (3.62) of [7]: if $f\in L^2\cap L^\infty$, then for all $s>0$
$$|D(P_{-s}f)|^2=\G(P_{-s} f,P_{-s} f)  \txt{$m$ a. e..}  \eqno (1.11)$$
Moreover, there is an increasing function  
$$\fun{\a}{(0,+\infty)}{(0,+\infty)}$$
such that, for all $u\in L^\infty(S,m)$ and $s>0$,
$$\a(s)||\G(P_{-s}u,P_{-s}u)||_{L^\infty}^\2\le ||u||_{L^\infty}   .  \eqno (1.12)$$
Another fact we need is the chain rule, or formula (2.12) of [7]: if $f,g\in\dc(\ec)$ and $\eta\in Lip(\R)$, then 
$$\G(\eta(f),g)=\eta^\prime(f)\G(f,g)  .  \eqno (1.13)$$
The last fact we need follows from formula (2.58) of [7]: if $f\in\dc(\ec)$, 
$$|DP_{-s}f|_w^2\le e^{-2K\b s}P_{-s}\G(f)  .  \eqno (1.14)$$

We briefly prove that, if $f\in\dc(\ec)$, then $P_{-\e}f\tends f$ in $\dc(\ec)$ as $\e\tends 0$. Since $\dc(\ec)$ is a Hilbert space, it suffices to prove that 

\noindent 1) $||P_{-\e}f||_{\dc(\ec)}\tends||f||_{\dc(\ec)}$ and

\noindent 2) $P_{-\e}f\tends f$ weakly in $\dc(\ec)$. 

By (1.14) we get that 
$$\limsup_{\e\searrow 0} ||P_{-\e}f||_{\dc(\ec)}\le ||f||_{\dc(\ec)}  .  $$
Thus, points 1) and 2) follow if we show that
$$\int_S\G(\phi,P_{-\e}f)\dr m\tends
\int_S\G(\phi,f)\dr m$$
for $\phi$ in a dense set of $\dc(\ec)$. The dense set is $V^2_\infty$; if $\phi\in V^2_\infty$, we get that $\D_\ec\phi\in L^\infty(S,m)\subset L^2(S,m)$ and the limit below follows; the equalities are (1.5). 
$$\int_S\G(\phi,P_{-\e}f)\dr m=
-\int_S\D_\ec\phi\cdot P_{-\e}f\dr m\tends
-\int_S\D_\ec\phi\cdot f\dr m=
\int_S\G(\phi,f)\dr m  .  $$

\vskip 2 pc
\centerline{\bf  \S 2}
\centerline{\bf Weak solutions of the continuity and Fokker-Planck equations}
\vskip 1 pc

A result of [3] says that a curve of measures $\mu_s$ in $\R^d$ is 2-absolutely continuous for the 2-Wasserstein distance if and only if it is a weak solution of the continuity equation for a drift $X\in L^2([t,0]\times\R^d,\L^1\otimes\mu_s)$. This result has been extended to Banach spaces in [14] and to measured metric spaces in [13] (see also sections 6 and 8 of  [5]). In proposition 2.4 below, we recall a version of theorem 3.5 of [13], which we shall state in the less general setting of $RCD(K,\infty)$ spaces. The advantage of this reduction in generality is that the key lemma 2.4 below holds also for curves of measures with unbounded densities; this allows us to extend [13] to unbounded densities, though in the much less general situation of $RCD(K,\infty)$ spaces. 

\noindent{\bf The test functions.} Let $t<0$; we consider the space $C^1([t,0],V^1_\infty)$ with the norm
$$||u||_{C^1([t,0],V^1_\infty)}=
\sup_{s\in[t,0]}||u_s||_{V^1_\infty}+
\sup_{s\in[t,0]}||\frac{\dr}{\dr s}u_s||_{V^1_\infty}    $$
and $C([t,0],V^2_\infty)$ with the norm
$$||u||_{C([t,0],V^2_\infty)}=\sup_{s\in[t,0]}||u_s||_{V^2_\infty}  .  $$
Note that, since $V^1_\infty$ and $V^2_\infty$ are Banach spaces, $C^1([t,0],V^1_\infty)$ and $C([t,0],V^2_\infty)$ are Banach spaces too. 

We shall say that $u$ is a test function, or that $u\in\tc$ for short, if 
$$u\in C^1([t,0],V^1_\infty)\cap C([t,0],V^2_\infty)  .  $$

By [7], the functions in $V^1_\infty$ are continuous; thus, if $u\in\tc$, then 
$u\in C([t,0],C(S))$ too.

\noindent{\bf The admissible curves.} Let $\fun{\mu}{[t,0]}{\pc(S)}$ be a function; we say that $\mu$ is admissible if 

\noindent 1) $\mu$ is continuous and 

\noindent 2) $\mu_s=\r_s m$ for all $s\in[t,0]$. 


\noindent{\bf Notation.} If $\fun{\mu}{[t,0]}{\pc(S)}$ is admissible, we shall always denote by $\r_s$ the density of $\mu_s$.

\noindent{\bf The space of drifts.} Throughout this section we follow [5], [11] and [13], which define the drift of a curve of measures $\fun{\mu}{[t,0]}{\pc(S)}$ as a linear operator on test functions, bounded for the seminorm 
$$\int_t^0\dr s\int_S\G(\phi_s,\phi_s)\dr\mu_s  .  $$
Since we work on $RCD(K,\infty)$ spaces, we shall follow the approach of [5] and identify  the drift operator with an element of a (rather abstract) Hilbert space which we call 
$\vc(\mu)$. 

Let $\phi,\psi\in\tc$; since $\G$ is bilinear and semi-positive-definite, Cauchy-Schwarz implies the inequality below. 
$$|\G(\phi_s,\psi_s)|\le\G(\phi_s,\phi_s)^\2\cdot\G(\psi_s,\psi_s)^\2  
\txt{$m$-a.e..}  \eqno (2.1)$$
As a consequence, if $\phi,\psi\in\tc$ we have that the function
$$\fun{A}{(s,x)}{\G(\phi_s,\psi_s)(x)}$$
belongs to $L^\infty([t,0]\times S,\L^1\otimes m)$. This implies that, if 
$\fun{\mu}{[t,0]}{\pc(S)}$ is admissible, the integral on the right in the formula below converges, allowing us to define the "inner product" $\inn{\phi}{\psi}_{\vc(\mu)}$. 
$$\inn{\phi}{\psi}_{\vc(\mu)}=
\int_t^0\dr s\int_S\G(\phi_s,\psi_s)\dr\mu_s  .  \eqno (2.2)$$
Though we called $\inn{\phi}{\psi}_{\vc(\mu)}$ an inner product, it does not separate points: a typical example is when there is $A\subset S$ with $m(A)>0$ but $\mu_s(A)=0$ for all 
$s\in[t,0]$. In this case, if $\phi_s$ is supported on $A$ for all $s\in[t,0]$, we get that 
$\inn{\phi}{\psi}_{\vc(\mu)}=0$ for all $\psi\in\tc$. Another example is the couple $\phi$, 
$\phi+c$ where $c$ is a constant: this couple is not separated by 
$\inn{\cdot}{\cdot}_{\vc(\mu)}$. 

Thus,  in order to have a Hilbert space, a further step is in order. First of all, we define the seminorm
$$||\phi||^2_{\vc(\mu)}\colon =\inn{\phi}{\phi}_{\vc(\mu)}    $$
and we say that $\phi\simeq\psi$ if 
$$||\phi-\psi||_{\vc(\mu)}=0  .   \eqno (2.3)$$
We shall call $\vc(\mu)$ the completion of $\frac{\tc}{\simeq}$ with respect to 
$||\cdot||_{\vc(\mu)}$. It is easy to see that the inner product of (2.2) depends only on the equivalence classes of $\phi$ and $\psi$ and separates points on $\frac{\tc}{\simeq}$; thus, $\vc(\mu)$ is a Hilbert space. 

We also note that (2.2) implies that, if $\phi\in\tc$, then 
$$||\phi||_{\vc(\mu)}\le
|t|^\2\cdot
||\phi||_{C([t,0],V^1_\infty)}  .  \eqno (2.4)$$

We can do the same construction for a measure $\eta\in\pc(S)$, provided it has a density. Indeed, let $\eta=\r m\in\pc(S)$ and let $\psi,\phi\in V^1_\infty$; we define
$$\inn{\phi}{\psi}_{\zc(\eta)}=
\int_S\G(\phi,\psi)\dr\eta=\int_S\G(\phi,\psi)\r\dr m    $$
which is well defined since $\G(\phi,\psi)\in L^\infty(S,m)$ by (2.1) and $\r\in L^1(S,m)$. 

We set $||\phi||^2_{\zc(\eta)}=\inn{\phi}{\phi}_{\zc(\mu)}$; given $\phi,\psi\in V^1_\infty$ we say that $\phi\simeq\psi$ if 
$$||\phi-\psi||_{\zc(\eta)}=0  .  $$
We shall call $\zc(\eta)$ the completion of $\frac{V^1_\infty}{\simeq}$ with respect to 
$||\cdot||_{\zc(\eta)}$; with the same argument as above, $\zc(\eta)$ is a Hilbert space. 

The following lemma relates $\vc(\mu)$ with $\zc(\mu_s)$. 

\lem{2.1} Let $\fun{\mu}{[t,0]}{\pc(S)}$ be admissible and let $V\in\vc(\mu)$. Then, for 
$\L^1$ a. e. $s\in[t,0]$, there is $V_s\in\zc(\mu_s)$ such that, for all $\phi\in\tc$, 
$$\inn{\phi}{V}_{\vc(\mu)}=
\int_t^0\inn{\phi_s}{V_s}_{\zc(\mu_s)}\dr s  \eqno (2.5)$$
and
$$||V||^2_{\vc(\mu)}=
\int_t^0||V_s||^2_{\zc(\mu_s)}\dr s  .  \eqno (2.6)  $$

\proof Let $\{ V_n \}_{n\ge 1}\subset\tc$ be a Cauchy sequence in $\vc(\mu)$; then,
$$\int_t^0\dr s\int_S\G(
V_n(s,\cdot)-V_m(s,\cdot),V_n(s,\cdot)-V_m(s,\cdot)
)(x)\dr\mu_s(x)  \tends 0$$
as $n,m\tends+\infty$. 

We find a subsequence $\{ V_{n_k} \}$ such that
$$||V_{n_{k+1}}-V_{n_k}||_{\vc(\mu)}\le\frac{1}{2^k}  . $$
Let us define a function $\fun{g}{[t,0]}{\R}$ by 
$$g(s)\colon=||V_{n_0}(s,\cdot)||_{\zc_s(\mu_s)}+
\sum_{k\ge 0}||V_{n_{k+1}}(s,\cdot)-V_{n_{k}}(s,\cdot)||_{\zc_s(\mu_s)}  \eqno (2.7)  $$
The formula below shows that $g$ belongs to $L^2(t,0)$; the first inequality is Minkowski's, the equality comes from the definition of $||\cdot||_{\vc(\mu)}$ and the second inequality from our choice of $V_{n_k}$.  
$$\left[
\int_t^0 g^2(s)\dr s
\right]^\2\le
\left[
\int_t^0||V_{n_0}(s,\cdot)||_{\zc_s(\mu_s)}^2\dr s
\right]^\2+
\sum_{k\ge 0}
\left[
\int_t^0||V_{n_{k+1}}(s,\cdot)-V_{n_k}(s,\cdot)||_{\zc_s(\mu_s)}^2\dr s
\right]^\2=$$
$$\left(
\int_t^0\dr s\int_S\G(V_{n_0}(s,\cdot),V_{n_0}(s,\cdot),)\dr\mu_s
\right)^\2+$$
$$\sum_{k\ge 0}\left(
\int_t^0\dr s
\int_S\G(V_{n_{k+1}}(s,\cdot)-V_{n_k}(s,\cdot),V_{n_{k+1}}(s,\cdot)-V_{n_k}(s,\cdot))
\dr\mu_s
\right)^\2\le$$
$$\left(
\int_t^0\dr s\int_S\G(V_{n_0}(s,\cdot),V_{n_0}(s,\cdot),)\dr\mu_s
\right)^\2+
\sum_{k\ge 1}\frac{1}{2^k}  <+\infty  .  $$
In particular, $g(s)$ is finite for a. e. $s\in(t,0)$. By (2.7) this implies that  
$\{ V_{n_k}(s,\cdot) \}$ is a Cauchy sequence in $\hat\zc_s(\mu_s)$ for a. e. $s\in(0,1)$; we call $V_s$ its limit. 

We prove that $V_s$ satisfies (2.5); we forgo the proof of (2.6) since it is similar. 
Since 
$$V_{n_k}=V_{n_0}+\sum_{j=1}^k[
V_{n_k}-V_{n_{k-1}}
]  ,  $$
formula (2.7) and the triangle inequality imply that 
$$||V_{n_k}(s,\cdot)||_{\zc_s(\mu_s)}\le g(s)
\qquad\forall k\ge 1\txt{a. e.} s\in(t,0)  .  \eqno(2.8)$$
Since $V_{n_k}(s,\cdot)\tends V_s$ in $\zc_s(\mu_s)$ for a. e. 
$s\in(t,0)$, we get that 
$$\int_S\G(\phi(s,\cdot),V_{n_k}(s,\cdot))\dr\mu_s\tends
\inn{\phi(s,\cdot)}{V_s}_{\zc_s(\mu_s)}  $$
for a. e. $s\in(t,0)$. By (2.8) and Cauchy-Schwarz we get that for a. e. $s\in(t,0)$ the first inequality below holds; the second inequality comes from the fact that $\phi\in\tc$ with a proof similar to (2.4). 
$$\left\vert
\int_S\G(\phi(s,\cdot),V_{n_k}(s,\cdot))\dr\mu_s
\right\vert  \le g(s)||\phi(s,\cdot)||_{\zc_s(\mu_s)}\le
Mg(s)  $$
By the last two formulas, dominated convergence and the fact that $g\in L^2(t,0)$, we get that 
$$\int_t^0\dr s\int_S\G(\phi(s,\cdot),V_{n_k}(s,\cdot))\dr\mu_s\tends
\int_t^0\inn{\phi(s,\cdot)}{V_s}_{\zc_s(\mu_s)}\dr s  .  $$
Since by assumption $V_{n_k}\tends V$ in $\vc(\mu)$, we get (2.5).

\fin

\vskip 1pc

\noindent{\bf The transport equations.}  Let $\fun{\mu}{[t,0]}{\pc(S)}$ be admissible and let 
$\fun{\L}{\tc}{\R}$ be a linear operator such that 
$$\L(\phi)=\L(\psi) \txt{if}\phi\simeq\psi    \eqno (2.9)$$
where $\simeq$ is the equivalence relation of (2.3). Following [11] and [13] we define
$$||\L||_{\vc^\prime(\mu)}=\sup\{
\L(\phi)\st \phi\in\tc \txt{and} ||\phi||_{\vc(\mu)}\le 1
\}   .   \eqno (2.10)$$
Clearly, if $||\L||_{\vc^\prime(\mu)}<+\infty$, then $\L$ can be extended to a bounded operator on the Hilbert space $\vc(\mu)$.  We also note that, if 
$||\L||_{\vc^\prime(\mu)}<+\infty$, then (2.9) holds. 

We say that the admissible $\fun{\mu}{[t,0]}{\pc(S)}$ is a weak solution of the continuity equation if the operator 
$$\hat\L_\mu(\phi)=
-\int_t^0\dr s\int_S
\partial_s\phi_s
\r_s\dr m  +
\int_S\phi_0\dr\mu_0-\int_S\phi_t\dr\mu_t     $$
belongs to $\vc^\prime(\mu)$. Note that all the integrals in the formula above converge: for instance, $\partial_s\phi_s$ is in $L^\infty$ by the definition of $\tc$, while $\mu_s=\r_sm$ is a probability measure.

By Riesz's representation theorem, $\mu$  is a weak solution of the continuity equation if and only if there is $Y\in\vc(\mu)$ such that for all 
$\phi\in\tc$ we have the first equality below; the second one is (2.5). 
$$\int_S\phi_0\dr\mu_0-\int_S\phi_t\dr\mu_t=
\int_t^0\dr s\int_S\partial_s\phi_s\dr\mu_s+\inn{Y}{\phi}_{\vc(\mu)}=
\int_t^0\dr s\int_S\partial_s\phi_s\dr\mu_s+
\int_t^0\inn{\phi_s}{Y_s}_{\zc(\mu_s)}\dr s  .  \eqno (2.11)$$

Let $\b>0$ be a diffusion coefficient that we fix once and for all; we say that $\mu$ is a weak solution of the Fokker-Planck equation if the operator   
$$\L_\mu(\phi)=
\int_t^0\dr s\int_S[
-\partial_s\phi_s-\frac{\b}{2}\D_\ec\phi_s
]\r_s\dr m  +
\int_S\phi_0\dr\mu_0-\int_S\phi_t\dr\mu_t     \eqno (2.12)$$
belongs to $\vc^\prime(\mu)$. Again, since $\phi\in\tc$ we can check that all the integrals in the formula above converge. 

Applying Riesz's representation theorem as above, $\mu$ is a weak solution of Fokker-Planck if and only if there is $V\in\vc(\mu)$ such that for all $\phi\in\tc$ we have the first equality below.
$$\int_S\phi_0\dr\mu_0-\int_S\phi_t\dr\mu_t=$$
$$\int_t^0\dr s\int_S\left[
\partial_s\phi_s+\frac{\b}{2}\D_\ec\phi_s
\right]\dr\mu_s+\inn{V}{\phi}_{\vc(\mu)}=$$
$$\int_t^0\dr s\int_S\left[
\partial_s\phi_s+\frac{\b}{2}\D_\ec\phi_s
\right]\dr\mu_s+
\int_t^0\inn{\phi_s}{V_s}_{\zc(\mu_s)}\dr s  .  \eqno (2.13)$$

Formulas (2.11) and (2.13) hold for the extrema $t$ and $0$; since admissible curves are continuous, they hold for all extrema $t\le s\le s^\prime\le 0$. We prove this for (2.13). Let 
$\fun{\eta^n}{\R}{\R}$ be a sequence of $C^\infty$ functions such that
$$0\le\eta^n(\tau)\le 1\qquad\forall\tau\in\R
\txt{and}
\eta^n(\tau)=0\txt{if}\tau\not\in\left[
s,s^\prime 
\right]  ,$$
$$\eta^n(\tau)\tends 1_{[s,s^\prime]}(\tau)\qquad\forall\tau\in\R  ,  $$
$$\dot\eta^n\L^1\tends\d_{s}-\d_{s^\prime}\txt{in}\pc(\R)    \eqno (2.14)$$
where $\d_s$ is the Dirac delta centred in $s$. 

If $\phi\in\tc$, then also $\eta^n\phi\in\tc$; taking it as a test function in (2.13) we have that
$$\int_S\eta^n_0\phi_0\dr\mu_0-\int_S\eta^n_t\phi_t\dr\mu_t=$$
$$\int_t^0\dr\tau\int_S[\dot\eta^n_\tau\phi_\tau+
\eta^n_\tau\partial_\tau\phi_\tau+
\frac{\b}{2}\eta^n_\tau\D_\ec\phi_\tau]\dr\mu_\tau+
\int_t^0\eta^n_\tau\inn{\phi_\tau}{V_\tau}_{\zc{(\mu_\tau)}}\dr\tau  .  \eqno (2.15)$$
We saw above that, if $\phi\in\tc$, then $\phi\in C([t,0],C(S))$; since $\mu$ is admissible we easily see that the map 
$$\fun{}{\tau}{
\int_S\phi_\tau\dr\mu_\tau
}   $$
is continuous; thus, by (2.14), 
$$\int_t^0\dot\eta^n_\tau\dr\tau\int_S\phi_\tau\dr\mu_\tau\tends
\int_S\phi_{s}\dr\mu_{s}-\int_S\phi_{s^\prime}\dr\mu_{s^\prime}  .  $$
On the other side, 
$$\int_t^0\eta^n_\tau\dr\tau
\int_S[\partial_\tau\phi_\tau+\inn{\phi_\tau}{V_\tau}_{Z(\mu_\tau)}+
\frac{\b}{2}\D_\ec\phi_\tau]\dr\mu_\tau
\tends
\int_s^{s^\prime}\dr\tau
\int_S[\partial_\tau\phi_\tau+\inn{\phi_\tau}{V_\tau}_{Z(\mu_\tau)}+
\frac{\b}{2}\D_\ec\phi_\tau ]\dr\mu_\tau     $$
by dominated convergence; to find the dominant it suffices to recall that, since $\phi\in\tc$, $\partial_\tau\phi_\tau$ and $\D_\ec\phi_\tau$ are bounded functions; moreover, since 
$V\in\vc(\mu)$ and $\phi\in\tc$ we have that 
$\inn{\phi_\tau}{V_\tau}_{\zc(\mu_\tau)}\in L^2(t,0)$ by Cauchy-Schwarz and (2.6). Lastly, 
$$\int_S\eta^n_0\phi_0\dr\mu_0-\int_S\eta^n_t\phi_t\dr\mu_t=0\qquad\forall n$$
since $\eta_0^n=\eta^n_t=0$. 

By (2.15) and the last three formulas we get that  
$$\int_S\phi_{s^\prime}\dr\mu_{s^\prime}-
\int_S\phi_{s}\dr\mu_{s}=
\int_s^{s^\prime}\dr\tau\int_S\left[
\partial_\tau\phi_\tau\dr\mu_\tau+\frac{\b}{2}\D_\ec\phi_\tau
\right]   +
\int_s^{s^\prime}\inn{\phi_\tau}{V_\tau}_{\zc(\mu_\tau)}\dr\tau  .  $$
Analogously, we have that, for all $t\le s<s^\prime\le 0$,
$$\int_S\phi_{s^\prime}\dr\mu_{s^\prime}-
\int_S\phi_{s}\dr\mu_{s}=
\int_s^{s^\prime}\dr\tau\int_S
\partial_\tau\phi_\tau\dr\mu_\tau   +
\int_s^{s^\prime}\inn{\phi_\tau}{V_\tau}_{\zc(\mu_\tau)}\dr\tau  .  \eqno (2.16)$$

\noindent{\bf Absolutely continuous curves of measures.} As in [3] and [13], the following fact is essential; it was proven in [14]. 

\lem{2.2} Let $\mu\in AC^2([0,1],\pc(S))$. Then, there is a Borel probability measure $\pi$ on $C([0,1],S)$ such that the three points below hold. 

\noindent $i$) $\pi$ concentrates on $AC^2([0,1],S)$;

\noindent $ii$) $\mu_t=(e_t)_\sharp\pi$. 

\noindent $iii$) Defining $|\dot\mu_t|$ as in section 1, 
$$\int_{C([0,1],S)}|\dot\g_t|^2\dr\pi(\g)=|\dot\mu_t|^2
\txt{for $\L^1$-a. e. $t\in[0,1]$.}  $$

\rm

\lem{2.3} Let $\mu\in AC^2([0,1],\pc(S))$ be admissible. Let $\phi\in V^1_\infty$, let 
$s\in[0,1)$, $h\in(0,1-s)$ and let $\pi$ be as in lemma 2.3 above. Then,
$$\left\vert
\int_{C([0,1],S)}[\phi(\g_{s+h})-\phi(\g_s)]\dr\pi(\g)
\right\vert\le
\int_{C([0,1],S)}\dr\pi(\g)\int_s^{s+h}|D\phi|_w(\g_\tau)\cdot|\dot\g_\tau|\dr\tau  .  
\eqno (2.17)$$

\proof As in formula (4.30) of [7], the fact that the elements of $V^1_\infty$ are Lipschitz immediately implies that 
$$\left\vert
\int_{C([t,0],S)}  [
\phi(\g_{s+h})-\phi(\g_s)
]  \dr\pi(\g)
\right\vert \le
\int_s^{s+h}\dr\tau\int_{C([t,0],S)}
|D\phi|(\g_\tau)\cdot ||\dot\g_\tau||\dr\pi(\g)    \eqno (2.18)$$
where $|D\phi|$ has been defined in (1.1). 

We must prove that we can substitute $|D\phi|_w(s,\g_\tau)$ to $|D\phi|(s,\g_\tau)$ in the formula above.  

We recall formula (1.11) which says that, for all $\e>0$ and for all 
$f\in L^\infty$,
$$|D{P_{-\e}f}|(x)=|D{P_{-\e}f}|_w(x)
\txt{for $m$-a. e. $x\in S$.}   $$
Since $\phi\in L^\infty$, we get that
$$|D{P_{-\e}\phi}|(x)=|D{P_{-\e}\phi}|_w(x)
\txt{for $m$-a. e. $x\in S$.}   $$
Applying (2.18) to $P_{-\e}\phi$ we get that, for $\e>0$, 
$$\left\vert
\int_{C([t,0],S)}  \{
[  P_{-\e}   
\phi
](\g_{s+h})-[
P_{-\e}\phi
](\g_s)
\}   \dr\pi(\g)
\right\vert \le
\int_s^{s+h}\dr\tau\int_{C([t,0],S)}
|D[P_{-\e}\phi]|_w(\g_\tau)\cdot ||\dot\g_\tau||\dr\pi(\g)  .  $$
Thus, (2.17) follows if we prove that, as $\e\tends 0$,
$$\left\vert
\int_{C([t,0],S)}  \{
[  P_{-\e}   
\phi
](\g_{s+h})-[
P_{-\e}\phi
](\g_s)
\}   \dr\pi(\g)
\right\vert  \tends
\left\vert
\int_{C([t,0],S)}  \{     
\phi (\g_{s+h})-\phi(\g_s)
\}   \dr\pi(\g) 
\right\vert   \eqno (2.19)$$
and 
$$\int_s^{s+h}\dr\tau
\int_{C([t,0],S)} \Big\vert
|D(P_{-\e}\phi)|_w(\g_\tau)-|D\phi|_w(\g_\tau)
\Big\vert    \cdot  ||\dot\g_\tau||\dr\pi(\g)\tends 0  \eqno (2.20)$$
As for (2.19), the first equality below comes from point $ii$) of lemma 2.2, i. e. that 
$(e_s)_\sharp\pi=\mu_s$; the second one comes from the fact that $P_{-\e}$ and $H_\e$ are in duality; the limit comes from the fact that $\phi$ is continuous and $H_\e\mu\tends\mu$ for $\e\tends 0$. 
$$\int_{C([t,0],S)}[P_{-\e}\phi](\g_s)\dr\pi(\g)=
\int_S[P_{-\e}(\phi)]\dr\mu_s=$$
$$\int_S\phi\dr(H_\e\mu_s)\tends 
\int_S\phi\dr\mu_s  \txt{as}\e\tends 0.  $$
Applying the same argument to $P_{-\e}\phi(\g_{s+h})$ we get (2.19). 

We prove (2.20). The inequality below comes from H\"older and points $ii$) and $iii$) of lemma 2.2. 
$$\int_s^{s+h}\dr\tau
\int_{C([t,0],S)} \Big\vert
|D(P_{-\e}\phi)|_w(\g_\tau)-|D\phi|_w(\g_\tau)
\Big\vert    \cdot  ||\dot\g_\tau||\dr\pi(\g)\le$$
$$\left[
\int_s^{s+h}\dr\tau\int_S   \Big\vert
|D(P_{-\e}\phi)|_w(x)-|D\phi|_w(x)
\Big\vert^2     \dr\mu_\tau(x)
\right]^\2  \cdot
\left[
\int_s^{s+h}||\dot\mu_\tau||^2\dr\tau
\right]^\2  .  $$
Since $\mu\in AC^2([t,0],\pc(S))$, we get that 
$$\int_s^{s+h}||\dot\mu_\tau||^2\dr\tau   <+\infty   .  $$
Thus, (2.20) follows if we prove that 
$$\int_s^{s+h}\dr\tau\int_S    \Big\vert
|D (P_{-\e}\phi)|_w(x)-|D \phi|_w(x)
\Big\vert^2   \dr\mu_\tau\tends
0  \txt{as $\e\tends 0$.}  \eqno (2.21)$$
We begin to note that
$$\int_S   \Big\vert
|D(P_{-\e}\phi)|_w(x)-|D\phi|_w(x)
\Big\vert^2\dr m(x)\tends 0
\txt{as}\e\tends 0    $$
because we saw at the end of section 1 that $P_{-\e}\phi\tends\phi$ in $\dc(\ec)$. Thus, for all sequences 
$\e_n\tends 0$, we can find a subsequence $\e_{n^\prime}$ such that
$$|D(P_{-\e_{n^\prime}}\phi)|_w(x)\tends|D\phi|_w(x)  
\txt{$m$ a. e. as $n^\prime\tends+\infty$.}  $$
Since $\mu$ is admissible, the convergence above is $\L^1\otimes\mu_s$-a. e.. Since 
$\phi\in V^1_\infty$, we can find $M>0$ such that the second inequality below holds, while the first one comes from (1.14).
$$|D(P_{-\e}\phi)|_w(x)\le 
e^{-2K\e}P_{-\e}(|D\phi|_w)(x)\le M\qquad\forall\e\in(0,1]  .  $$
Now (2.21) follows from the last two formulas and dominated convergence for the measure $\L^1\otimes\mu_s$. 

\fin

The next proposition follows from lemmas 2.2 and 2.3 exactly as in [3] and theorem 3.5 of [13]; thus, we forego its proof.  

\prop{2.4} Let $\fun{\mu}{[t,0]}{\pc(S)}$ be admissible. Then, $\mu\in AC^2([t,0],\pc(S))$ if and only if it is a weak solution of the continuity equation with drift $Y\in\vc(\mu)$. 

Moreover, we have that
$$\int_t^0 ||\dot\mu_s||^2\dr s= ||Y||_{\vc(\mu)}^2  .  $$

\rm

\vskip 2pc
\centerline{\bf \S 3}
\centerline{\bf Smooth approximation.}
\vskip 1pc

In this section we approximate the solutions of the Fokker-Planck equation with curves of measures having smoother densities and drifts. Our method is the standard one (see for instance [1]).

Let $\fun{\mu}{[t,0]}{\pc(S)}$ be admissible and let us suppose that (5) holds. 

\noindent{\bf Step 1. Regularisation in space.} In this step, we approximate $\mu_s$ with a curve of measures having a density $\hat\r_s^\e\in V^2_\infty$, but only Borel regular in $s$. 

Following [7], we choose $k\in C^\infty_c(0,+\infty)$ with $k\ge 0$ and
$$\int_0^{+\infty}k(r)\dr r=1  .  \eqno (3.1)$$
Given $\phi\in L^2(S,m)$, for $\e>0$ we set 
$$M^\e\phi=\frac{1}{\e}
\int_0^{+\infty}(P_{-r}\phi)k\left(\frac{r}{\e}\right)\dr r  .  \eqno (3.2)$$
By [7], we have that $M^\e\phi\in\dc(\Delta_\ec)$ and 
$$\frac{\b}{2}\Delta_\ec (M^\e\phi)=\frac{-1}{\e^2}
\int_0^{+\infty}(P_{-r}\phi)k^\prime\left(\frac{r}{\e}\right)\dr r    \eqno (3.3)$$
(we recall that $P_{-r}$ is the flow of $\frac{\b}{4}\ec$.)

\lem{3.1} For all $\e>0$ there is a constant $D_1(\e)>0$ such that the following holds. Let   
$\phi\in L^\infty(S,m)$ and let 
$$\hat\phi^\e\colon= M^\e\phi   $$ 
with $M^\e$ defined as in (3.2); then
$$||\G(M^\e\phi,M^\e\phi)||_{\infty}\le 
D_1(\e)||\phi||_\infty  \eqno (3.4)$$
and
$$||\D_\ec (M^\e\phi)||_{L^\infty}\le D_1(\e)||\phi||_\infty   .   \eqno (3.5)$$
Moreover, the operator $M^\e$ is self-adjoint: if $\phi,\psi\in L^2(S,m)$, then 
$$\int_S\psi\cdot M^\e\phi\dr m=
\int_S M^\e\psi\cdot\phi\dr m  .  \eqno (3.6)$$
We also have 
$$\D_\ec(M^\e\phi)=M^\e(\D_\ec \phi)   \qquad
\forall\phi\in\dc(\D_\ec)  .  \eqno (3.7)$$

\proof Since $k\in C_c(0,+\infty)$, for some $\d>0$ it is supported in $[\d,\frac{1}{\d}]$; thus, $k\left(\frac{r}{\e}\right)$ is supported in $\left[ \e\d,\frac{\e}{\d} \right]$ and (3.2) becomes
$$M^\e\phi=\frac{1}{\e}
\int_{\d\e}^{\frac{\e}{\d}}(P_{-r}\phi)k\left(\frac{r}{\e}\right)\dr r  .  \eqno (3.8)$$
It is easy to see that $||\G(\phi,\phi)||_\infty^\2$ is convex and lower semicontinuous for the topology of $\dc(\ec)$; we briefly prove how this implies that it satisfies Jensen's inequality. We fix $\tilde\phi\in\dc(\ec)$ and $\l>0$; since $\dc(\ec)$ with the inner product 
$\inn{\cdot}{\cdot}_{\dc(\ec)}$ of section 1 is Hilbert and  the epigraph of 
$||\G(\cdot,\cdot)||_\infty^\2$ is closed, we can find $v_\l\in\dc(\ec)$ which separates the epigraph from $(\tilde\phi,||\G(\tilde\phi,\tilde\phi||_\infty^\2-\l))$ in $\dc(\ec)\times\R$. In other words, we have that, for all $\psi\in\dc(\ec)$,
$$||\G(\psi,\psi)||_\infty^\2-
||\G(\tilde\phi,\tilde\phi)||_\infty^\2-
\inn{\tilde\psi-\phi}{v_\l}_{\dc(\ec)}\ge -\l  .  $$
Now we take $\tilde\phi=M_\e\phi$, $\psi=P_{-r}\phi$; we multiply the formula above by 
$\frac{1}{\e}k\left(\frac{r}{\e}\right)$, integrate on $(0,+\infty)$ and recall (3.1), getting the inequality below. 
$$||\G(M_\e\phi,M_\e\phi)||^\2_\infty\le
\frac{1}{\e}\int_{\e\d}^\frac{\e}{\d}\left\vert\left\vert
\G(P_{-r},P_{-r})
\right\vert\right\vert_\infty^\2    k\left(\frac{r}{\e}\right)\dr r
    +\l  .  $$
If we let $\l\tends 0$, we get the first inequality below; the second one is (1.12). 
$$||\G(M^\e\phi,M^\e\phi)||_\infty^\2\le
\sup_{r\in [{\d\e},\frac{\e}{\d}]}
||\G(P_{-r}\phi,P_{-r}\phi)||_\infty^\2 \le
\frac{1}{a(\e\d)}||\phi||_\infty  .  $$
This proves (3.4). Next, we recall that, by the maximum principle, 
$$||P_{-r}\phi||_\infty\le||\phi||_\infty  .  $$
Together with (3.3), this implies (3.5). 

Formula (3.6) follows easily from Fubini's theorem and the fact that $P_{-r}$ is self-adjoint. 

Next, we recall that $P_{-r}$ is the semigroup associated to the maximal monotone operator $\D_\ec$; in particular ([9]), if $\phi\in\dc(\D_\ec)$ and $r>0$ we have 
$\D_\ec P_{-r}\phi=P_{-r}\D_\ec\phi$. Together with (3.8) this easily implies (3.7). 

\fin

Let $\mu_s=\r_s m$ be a continuous curve of measures; we set 
$\hat\r^\e_s\colon= M^\e\r_s$ and $\hat\mu^\e_s\colon=\hat\r^\e_s m$. By (3.4) and (3.5), 
$\hat\r^\e_s$ is a bounded, Borel function from $[t,0]$ to $V^2_\infty$. 

\noindent{\bf Step 2. Definition of the drift.} Let us suppose that $\mu$ solves the Fokker-Planck equation for a drift $V\in\vc(\mu)$ and the continuity equation for a drift 
$Y\in\vc(\mu)$; let the curve $\hat\mu^\e$ be defined as at the end of the last step. We want to find the drifts, which we shall call $\hat V^\e$ and $\hat Y^\e$ respectively, of the Fokker-Planck and continuity equations satisfied by 
$\hat\mu^\e$.

In order to find the drifts of $\hat\mu^\e$, for $\phi\in L^2(S,m)$ we set  
$$\hat V_s^\e(\phi)=
\inn{ V_s}{
M^\e\phi
}_{\zc(\mu_s)}$$
and 
$$\hat Y^\e(\phi)=
\inn{ Y_s}{
M^\e\phi
}_{\zc(\mu_s)}   .   $$
Let $\phi\in\tc$; the first equality below is the definition of $\hat V^\e$; the first inequality is Cauchy-Schwarz; the second one follows from (3.8) and the fact that $\G$ is convex; it is the same argument we used in lemma 3.1. The third inequality is (1.14); the second equality comes from Fubini, the definition of $\hat\r_s$ and the fact that $P_{-r}$ is self-adjoint; the last equality is the definition of the norm in 
$\zc(\hat\mu_s^\e)$. 
$$\inn{\hat V^\e_s}{\phi_s}_{\zc(\hat\mu_s)}=
\inn{V_s}{M^\e\phi_s}_{\zc(\mu_s)}\le$$
$$||V_s||_{\zc(\mu_s)}\cdot\left[
\int_S\G(M^\e\phi_s,M^\e\phi_s)\r_s\dr m
\right]^\2  \le$$
$$||V_s||_{\zc(\mu_s)}\cdot\left[
\int_S\r_s\dr m\frac{1}{\e}
\int_{\e\d}^\frac{\e}{\d}\G(P_{-r}\phi_s,P_{-r}\phi_s)k\left(\frac{r}{\e}\right)\dr r
\right]^\2  \le$$
$$||V_s||_{\zc(\mu_s)}\cdot e^{-K\b\e}\cdot\left[
\int_S\r_s\dr m\frac{1}{\e}
\int_{\e\d}^\frac{\e}{\d}P_{-r}\G(\phi_s,\phi_s)k\left(\frac{r}{\e}\right)\dr r
\right]^\2  =$$
$$||V_s||_{\zc(\mu_s)}\cdot e^{-K\b\e}\cdot\left[
\int_S\hat\r^\e_s\G(\phi_s,\phi_s)\dr m
\right]^\2=$$
$$||V_s||_{\zc(\mu_s)}\cdot e^{-K\b\e}\cdot
||\phi||_{\zc(\hat\mu^\e)}  .  $$
From this and (2.6) we get that 
$$||\hat V||_{\vc(\hat\mu^\e)}\le 
e^{-K\b\e}\cdot ||V||_{\vc(\mu)}   .  \eqno (3.9)$$
Analogously, we get that 
$$||\hat Y^\e||_{\vc(\hat\mu^\e)}\le
e^{-K\b\e}||Y||_{\vc(\mu)}  .  \eqno (3.10)$$

\noindent{\bf Step 3. The equation.} In this step we show that $\hat\mu^\e$ satisfies the Fokker-Planck equation with drift $\hat V^\e$ and the continuity equation with drift 
$\hat Y^\e$. Let $\phi\in\tc$; the first equality below comes from the definition of 
$\hat\mu^\e_s$ and of the drift $\hat V^\e$; the second one comes from (3.6). For the third one, we use (3.7). The fourth and last equality comes from the fact that $\mu$ satisfies the Fokker-Planck equation with drift $V$, and $M^\e\phi\in\tc$.  
$$\int_t^0\dr s\int_S[
\partial_s\phi_s+\frac{\b}{2}\D_\ec\phi_s
]  \dr \hat\mu_s+
\inn{\hat V^\e}{\phi}_{\vc(\hat\mu^\e)}-
\int_S\phi_0\dr \hat\mu_0+\int_S\phi_t\dr \hat\mu_t=$$
$$\int_t^0\dr s\int_S[
\partial_s\phi_s+\frac{\b}{2}\D_\ec\phi_s
] 
\cdot M^\e\r_s\dr m  +
\inn{ V}{
M^\e\phi
}_{\vc(\mu)}-$$
$$\int_S\phi_0\cdot M^\e\r_0\dr m+
\int_S\phi_t\cdot M^\e\r_t\dr m=$$
$$\int_t^0\dr s\int_S[
\partial_s(M^\e\phi_s)+
\frac{\b}{2}M^\e(\D_\ec\phi_s)
]  \r_s\dr m+
\inn{V}{M^\e\phi}_{\vc(\mu)}-$$
$$\int_S(M^\e\phi_0) \cdot \r_0\dr m+\int_SM^\e(\phi_t)\r_t\dr m  =$$
$$\int_t^0\dr s\int_S[
\partial_s(M^\e\phi_s)+
\frac{\b}{2}\D_\ec(M^\e\phi_s)
]  \r_s\dr m+
\inn{V}{M^\e\phi}_{\vc(\mu)}-$$
$$\int_S(M^\e\phi_0)\cdot\r_0\dr m+\int_S(M^\e\phi_t)\cdot\r_t\dr m=0  .  \eqno (3.11)$$
The proof for the continuity equation is similar. 

\noindent{\bf Step 4. Regularisation in time.} In this step, we regularise in time the density 
$\hat\r^\e$; at the end, we shall get a density $\r^\e\in\tc$. The reason we do this is that, as we shall see in lemma 3.3 below, the curve of measures $\fun{}{s}{\r^\e_s m}$ is Lipschitz for the 2-Wasserstein distance; this is one of the hypotheses of the integration by parts formula in lemma 4.2. 

Let $\eta\in C^\infty(\R,\R)$ be supported in 
$(-1,1)$; let it be even, nonnegative and with integral 1. We set
$$\eta^\e(r)=\frac{1}{\e}\eta(\frac{1}{\e} r)  .   $$
We set $\hat\r^\e_\tau=\hat\r^\e_t$ if $\tau\le t$ and $\hat\r^\e_\tau=\hat\r^\e_0$ if 
$\tau\ge 0$; we define $\mu^\e_s$ and $\r^\e_s$ as 
$$\mu^\e_s=\r^\e_s m=
\left[
\int_\R\eta^\e(\tau-s)\hat\r^\e_\tau\dr\tau
\right]  m .   \eqno (3.12)$$
Again setting 
$\hat V^\e_\tau=0=\hat Y^\e_\tau$ if $\tau\not\in[t,0]$, we define 
$V^\e_s,Y^\e_s\in\zc(\mu^\e_s)$ as
$$V^\e_s=\int_\R\eta^\e(\tau-s)\hat V^\e_\tau\dr\tau,\qquad
Y^\e_s=\int_\R\eta^\e(\tau-s)\hat Y^\e_\tau\dr\tau  .  $$
In other words, if $\phi\in V^1_\infty$,
$$\inn{V^\e_s}{\phi}_{\zc(\mu^\e_s)}=
\int_\R\eta^\e(\tau-s)\inn{\hat V^\e_\tau}{\phi}_{\zc(\hat\mu^\e_\tau)}\dr\tau ,\qquad
\inn{Y^\e_s}{\phi}_{\zc(\mu^\e_s)}=
\int_\R\eta^\e(\tau-s)\inn{\hat Y^\e_\tau}{\phi}_{\zc(\hat\mu^\e_\tau)}\dr\tau .  
\eqno (3.13)$$
We want to estimate the norm of $V^\e_s$; the first inequality below comes from (3.13) and Cauchy-Schwarz on $\zc(\mu^\e_\tau)$, the second one is Cauchy-Schwarz on $L^2(\R)$; the last equality is the definition of $\mu^\e_s$. 
$$\inn{V^\e_s}{\phi}_{\zc(\mu^\e_s)}\le
\int_\R\eta^\e(\tau-s)||\hat V^\e_\tau||_{\zc(\hat\mu^\e_\tau)}
\left(
\int_S\G(\phi,\phi)\dr\hat\mu^\e_\tau
\right)^\2   \dr\tau  \le$$
$$\left[
\int_\R\eta^\e(\tau-s)||\hat V^\e_\tau||_{\zc(\hat\mu^\e_\tau)}^2\dr\tau
\right]^\2\cdot
\left[
\int_\R\eta^\e(\tau-s)\dr\tau\int_S\G(\phi,\phi)\dr\hat\mu^\e_\tau
\right]^\2=$$
$$\left[
\int_\R\eta^\e(\tau-s)||\hat V^\e_\tau||_{\zc(\hat\mu^\e_\tau)}^2\dr\tau
\right]^\2 \cdot
\left[
\int_\R\G(\phi,\phi)\dr\mu^\e_s
\right]^\2  .  $$
By the definition of the norm in $\zc(\mu^\e_s)$ this implies that 
$$||V^\e_s||_{\zc(\mu^\e_s)}\le
\left[
\int_\R\eta^\e(\tau-s)||\hat V^\e_\tau||_{\zc(\hat\mu^\e_\tau)}^2\dr\tau
\right]^\2    .  \eqno (3.14)$$
The first equality below is (2.6), the first inequality comes from the formula above; the second equality comes from Fubini and the fact that $\eta^\e$ has integral 1; the second inequality is (3.9). 
$$||V^\e||_{\vc(\mu^\e)}^2=
\int_t^0||V^\e_s||^2_{\zc(\mu^\e_s)}\dr s\le
\int_t^0\dr s\int_\R\eta^\e(\tau-s)||\hat V^\e_\tau||_{\zc(\hat\mu^\e_\tau)}^2\dr\tau=$$
$$||\hat V^\e||^2_{\vc(\hat\mu^\e)}\le
e^{-K\b\e}||V||^2_{\vc(\mu)}  .  \eqno (3.15)$$
Analogously, 
$$|| Y^\e ||_{\vc(\mu^\e)}^2\le
e^{-K\b\e}||Y||^2_{\vc(\mu)}  . $$
Since $||\hat V^\e_\tau||_{\zc(\hat\mu^\e_\tau)}^2\in L^1(t,0)$ by (2.6) and 
$\eta_\e\le\frac{D_1}{\e}$ by definition, (3.14) implies the first inequality below, while the second one comes from (3.9). 
$$||V^\e_s||_{\zc(\mu^\e_s)}\le
\sqrt{\frac{D_1}{\e}}||\hat V^\e||_{\vc(\hat\mu^\e)}\le
\sqrt\frac{D_1}{\e}\cdot e^{-K\b\e}\cdot||V||_{\vc(\mu)}  .  \eqno (3.16)$$
Analogously, 
$$||Y^\e_s||_{\zc(\mu^\e_s)}\le
\sqrt{\frac{D_1}{\e}}\cdot e^{-K\b\e}\cdot||Y||_{\vc(\mu)}   .    $$

We saw in step 3 above that $\hat\mu^\e$ satisfies the Fokker-Planck equation with drift 
$\hat V^\e$; by the way we extended $\hat\mu^\e$ and $\hat V^\e$ outside $[t,0]$, we see that $(\hat\mu^\e,\hat V^\e)$ satisfies Fokker-Planck over all $\R$. Integrating against 
$\eta^\e$ we get as in (3.11) that
$$\int_t^0\dr s\int_S[
\partial_s\phi_s+\frac{\b}{2}\D_\ec\phi_s
]\r^\e_s\dr m+
\inn{V^\e}{\phi}_{\vc(\mu)}  =   
\int_S\phi_0\dr\mu^\e_0-\int_S\phi_t\dr\mu^\e_t  .  \eqno (3.17)$$
Similarly, $\mu^\e$ satisfies the continuity equation with drift $Y^\e$.

This discussion brings us to the following lemma. 

\lem{3.3} Let the admissible $\fun{\mu}{[t,0]}{\pc(S)}$ be a weak solution of the Fokker-Planck equation for a drift $V\in\vc(\mu)$. For $\e\in(0,1]$ let $\mu^\e=\r^\e m$ be defined as in (3.12) and let $V^\e\in\vc(\m^\e)$ be as in (3.13); let us suppose that (5) holds. Then, the following points hold.

\noindent 1) If $\e\in(0,1]$, then 
$$\frac{1}{C}\le\r_s^\e\le C  .  $$

\noindent 2) $\r^\e\in\tc$. 

\noindent 3) $\mu^\e$ satisfies the Fokker-Planck equation with drift $V^\e$. If moreover 
$\mu$ satisfies the continuity equation with drift $Y\in\vc(\mu)$, then $\mu^\e$ satisfies the continuity equation for the drift $Y^\e$ defined in (3.13).  

\noindent 4) $||V^\e||_{\vc(\mu^\e)}\le e^{-K\b\e}||V||_{\vc(\mu)}$.  If moreover 
$\mu$ satisfies the continuity equation with drift $Y\in\vc(\mu)$, we also have that
$||Y^\e||_{\vc(\mu^\e)}\le e^{-K\b\e}||Y||_{\vc(\mu)}$. 

\noindent 5) For all $s\in[t,0]$ we have 
$$|| V^\e_s ||_{\zc(\mu^\e_s)}\le
\sqrt{\frac{D_1}{\e}}\cdot e^{-K\b\e}||V||_{\vc(\mu)}     .   $$ 
If $\mu$ solves the continuity equation with drift $Y$, for all $s\in[t,0]$ we have
$$|| Y^\e_s ||_{\zc(\mu^\e_s)}\le
\sqrt{\frac{D_1}{\e}}\cdot e^{-K\b\e}||Y||_{\vc(\mu)}     .    $$

\noindent 6) $\r^\e\tends\r$ in $L^2([t,0]\times S,\L^1\otimes m)$  and 
$\mu^\e\tends\mu$ in $C([t,0],\pc(S))$ as $\e\tends 0$. 

\noindent 7) There is $\d(\e)\tends 0$ as $\e\tends 0$ such that, for all $s\in[t,0]$, 
$$\ent_m\m_s   \le
\ent_m\m^\e_s+\d(\e)  .  $$


\rm

\proof We sketch the proof of the inequality on the right of point 1); the inequality on the left is similar. The equality below is the definition of $\r^\e_s$ in (3.12); the first inequality follows since $\eta^\e$ is a probability density, the second one from (3.8) and the maximum principle for $P_{-r}$ and the third one from (5). 
$$|\r^\e_s(x)|=
\left\vert
\int_{\R}\eta^\e(\tau-s)\hat\r^\e_\tau(x)\dr\tau
\right\vert    \le$$
$$\sup_{\tau\in[t,0]}||\hat\r^\e_\tau||_\infty\le
\sup_{s\in[t,0]}||\r_s||_\infty\le C  .  $$

We prove point 2). From point 1), (3.4) and (3.5) we see that  the map 
$\fun{}{\tau}{\hat\r^\e_\tau}$ is Borel and bounded from $[t,0]$ to $V^2_\infty$. Since 
$\eta^\e\in C^\infty_0(\R)$, differentiation under the integral sign shows that 
$\fun{}{\tau}{\r_\tau^\e}$ is $C^1$ from $[t,0]$ to $V^2_\infty$; in particular, it belongs to 
$\tc$. 

Point 3) for the Fokker-Planck equation is formula (3.17); the continuity equation follows by a similar argument.  

Point 4) is (3.15) and the formula that follows it. Analogously, point 5) is formula (3.16) and the one that follows it. 

We prove point 6). Since $\r_s$ is bounded by (5) and $m$ is a probability measure, we have that $\r_s\in L^2(S,m)$; the strong continuity of $P_{-r}$ implies that, for all $s\in[t,0]$, 
$$\int_S|P_{-\e}\r_s-\r_s|^2\dr m\tends 0
\txt{as}
\e\tends 0  \qquad\forall s\in(t,0).  $$
Together with (3.8) this implies that 
$$\int_S|M^{\e}\r_s-\r_s|^2\dr m\tends 0
\txt{as}
\e\tends 0  \qquad\forall s\in(-t,0).  $$
By point 1) we have that 
$$\int_S|M^{\e}\r_s-\r_s|^2\dr m\le 4C^2  \qquad\forall s\in(t,0)$$
and thus dominated convergence implies convergence in 
$L^2([t,0]\times S,\L^1\otimes m)$. By convolution with $\eta^\e$ this implies that 
$\r^\e\tends\r$ in $L^2([t,0]\times S,L^1\otimes m)$. 

We prove convergence in $C([t,0],\pc(S))$. We begin to show that 
$$\hat\mu^\e_s=
\frac{1}{\e}\int_{\e\d}^\frac{\e}{\d}P_{-r}\r_sk\left( \frac{r}{\e} \right)\dr r \cdot m=
\frac{1}{\e}\int_{\e\d}^\frac{\e}{\d}(H_r\mu_s)\dr r
\tends\mu_s$$
uniformly in $s\in[-t,0]$ as $\e\tends 0$. 

First we note that, pointwise, $\hat\mu^\e_s\tends\mu_s$ for all $s$ by [7]; we omit the easy proof, based on the fact, which we saw after  formula (1.4), that the gradient flow 
$\fun{}{\e}{H_\e\mu_s}$ is continuous and exists for all initial conditions $\mu_s$.  Convergence is uniform by Ascoli-Arzel\`a if we prove that the maps 
$\fun{}{s}{\hat\mu_s}$ from $[t,0]$ to $\pc(S)$ have a modulus of continuity independent of $\e\in(0,1]$. We prove this: since $\mu$ is a continuous curve we can find a modulus of continuity $\omega$ for $\mu$. We recall that the map $\fun{}{\mu}{H_\e\mu}$ is Lipschitz with Lipschitz constant $e^{-\frac{K\e}{2}}$ by (1.4), implying that $H_\e\mu_s$ has modulus of continuity $\max(e^{-\frac{K\e}{2\d}},e^{-\frac{K\e\d}{2}})\cdot\omega$. For the convergence of $\mu^\e$, we convolute $\hat\mu^\e$ with $\eta^\e$ and use the convexity of the square of the Wasserstein distance; we easily see that $\mu^\e\tends\mu$ uniformly in $\pc(S)$. 

As for point 7), it follows immediately from point 6) and the fact that the entropy is lower semicontinuous for the Wasserstein distance. 


\fin

\vskip 2 pc
\centerline{\bf  \S 4}
\centerline{\bf Integration by parts}
\vskip 1 pc

In this section, we prove that (6) holds if the curve $\mu$ is "regular" in the sense of [5]. The heavy hauling will be done by point $iii$) of lemma 12.4 of [5], which we recall in lemma 4.2 below for convenience. 

\lem{4.1} Let $\mu\in AC^2([t,0],\pc(S))$ be admissible and let (5) hold; let us suppose that 
$$\int_t^0\dr s\int_S\frac{\G(\r_s,\r_s)}{\r_s}\dr m<+\infty  .  \eqno (4.1)$$
Then, the map $\fun{\Phi_\r}{\tc}{\R}$ defined by 
$$\fun{\Phi_\r}{\phi}{
\int_t^0\dr s\int_S\G(\r_s,\phi_s)\dr m  
}     $$
is bounded for the $||\cdot||_{\vc(\mu)}$ norm. As a consequence, $\Phi_\r$ can be extended to an element of $\vc^\prime(\mu)$; since $\vc(\mu)$ is a Hilbert space, there is 
$O(\r)\in\vc(\mu)$ such that 
$$\Phi_\r(V)=\inn{O(\r)}{V}_{\vc(\mu)}\qquad
\forall V\in\vc(\mu)  .  $$

\proof We prove that the map 
$$\fun{}{\phi}{
\int_t^0\dr s\int_S\G(\r_s,\phi_s)\dr m  
}     $$
is bounded. This follows by the inequality below (which is Cauchy-Schwarz) and (4.1).
$$\left\vert
\int_t^0\dr s\int_S\G(\r_s,\phi_s)\dr m
\right\vert   \le
\left[
\int_t^0\dr s\int_S\G(\log\r_s,\log\r_s)\r_s\dr m
\right]^\2\cdot
\left[
\int_t^0\dr s\int_S\G(\phi_s,\phi_s)\r_s\dr m
\right]^\2   .  $$

\fin

\lem{4.2} Let $\mu\in AC^2([t,0],\pc(S))$ be admissible and let us call $Y$ the drift of the continuity equation satisfied by $\mu$. Let us suppose that $\mu$ is regular in the sense of section 12 of [5], i. e. that $\fun{\mu}{[t,0]}{\pc(S)}$ is Lipschitz, the right hand side of (5) holds and there is $M>0$ such that 
$$\int_S\frac{\G(\r_s,\r_s)}{\r_s}\dr m\le M\qquad\forall s\in[t,0]  .  \eqno (4.2)$$

Then, we have that  
$$\ent_m(\mu_0)-\ent_m(\mu_t)=
\Phi_\r(Y)=
\int_t^0\inn{O_\tau(\r)}{Y_\tau}_{\zc(\mu_\tau)}\dr\tau   \eqno (4.3)$$
where the second equality comes from lemma 4.1. 

\rm

\lem{4.3} Let $\mu\in AC^2([t,0],\pc(S))$ be admissible; let $\fun{\mu}{[t,0]}{\pc(S)}$ be Lipschitz, let (5) and (4.2) hold. Let us suppose moreover that $\mu$ solves the Fokker-Planck equation for a drift $V\in\vc(\mu)$ and the continuity equation for a drift $Y\in\vc(\mu)$. 

Let the operator $\L_\mu$ be defined as in (2.12); then we have that
$$||\L_\mu||_{\vc^\prime(\mu)}^2=$$
$$\inn{Y}{Y}_{\vc(\mu)}   +   \eqno (4.4)_a$$
$$\b\ent_m(\mu_0)-\b\ent_m(\mu_t)+   \eqno (4.4)_b$$
$$\frac{\b^2}{4}\inn{O(\r)}{O(\r)}_{\vc(\mu)}  .  \eqno (4.4)_c$$

\vskip 1pc
\rm

\noindent{\bf Remark.} In [16], $Y$ is called the current velocity; its norm equals the metric velocity of the curve $\mu_s$ by proposition 2.4.
Formula $(4.4)_b$ is the contribution of the entropy; in $(4.4)_c$ we see the part of the kinetic energy due to $\frac{\nabla\r_s}{\r_s}$; in [16] this term is called the osmotic velocity, since it is the component of the velocity due to diffusion. Heuristically, the drift of Fokker-Planck is the current velocity plus the osmotic velocity and the formula above is simply the square of the binomial; $(4.4)_b$ represents the double product by lemma 4.2. 

\proof By point 1) of the hypotheses, $\mu$ satisfies the continuity equation (2.11) for a drift $Y\in\vc(\mu)$; from (2.11) and (2.12) we get that, for $\phi\in\tc$, 
$$\L_\mu(\phi)=
-\frac{\b}{2}\int_t^0\dr s\int_S\D_\ec\phi_s\dr\mu_s+
\inn{Y}{\phi}_{\vc(\mu)}  .  $$
Using (2.5), (2.10) and the integration by parts formula (1.5) we get that
$$||\L_\mu||_{\vc^\prime(\mu)}=$$
$$\sup\left\{
\frac{\b}{2}\int_t^0\dr s\int_S\G(\r_s,\phi_s)\dr m+
\int_t^0\inn{Y_{s}}{\phi_s}_{\zc(\mu_s)}\dr s\st
\phi\in\tc,\quad ||\phi||_{\vc(\mu)}\le 1
\right\}    .   \eqno(4.5)$$

We saw in proposition 2.4 that $Y\in\vc(\mu)$; together with (4.2), lemma 4.1 and the fact that $\tc$ is dense in $\vc(\mu)$ by definition, the last formula implies that
$$||\L_\mu||_{\vc^\prime(\mu)}^2=
\inn{Y+\frac{\b}{2}O(\r)}{Y+\frac{\b}{2}O(\r)}_{\vc(\mu)}   .  $$
By the properties of the internal product, we get that
$$||\L_\mu||_{\vc^\prime(\mu)}^2=||Y||_{\vc(\mu)}^2+
\b\inn{Y}{O(\r)}_{\vc(\mu)}+
\frac{\b^2}{4}\inn{O(\r)}{O(\r)}_{\vc(\mu)}  .  $$
Now (4.4) follows from lemma 4.2. 

\fin

\vskip 2 pc
\centerline{\bf  \S 5}
\centerline{\bf Proof of theorem 1}
\vskip 1 pc

In this section, we want to prove formula (6) with the hypothesis that $\mu$ satisfies the Fokker-Planck equation, but without knowing a priori that $\mu$ is regular. We begin recalling a standard lower semicontinuity result. 

\lem{5.1} Let $t<0$ and let $\fun{\mu^n,\mu}{[t,0]}{\pc(S)}$ be admissible curves. We suppose the following.

\noindent 1) For all $n\in\N$, $\mu^n$ satisfies the continuity equation and the Fokker-Planck equations with drifts respectively 
$$Y^n\in\vc(\mu^n)\txt{and}
V^n\in\vc(\mu^n) . $$

\noindent 2) $\mu^n\tends\mu$ in $C([t,0],\pc(S))$. 

\noindent 3) There is $D_3>0$ such that
$$||Y^n||_{\vc(\mu^n)}\le D_3\qquad\forall n\ge 1  .  $$

Then, $\mu$ is admissible, $\mu\in AC^2([t,0],\pc(S))$ and the following formula hold.
$$\int_t^0||\dot\mu_s||^2\dr s\le
\liminf_{n\tends+\infty}\int_t^0||\dot\mu^n_s||^2\dr s  .  \eqno (5.1)$$
If moreover $\mu^n=\r^n_sm$ with $\{\r^n_s\st s\in[t,0],n\ge 1\}$ uniformly integrable, then
$$||\L_\mu||_{\vc^\prime(\mu)}\le
\liminf_{n\tends+\infty}||\L_{\mu^n}||_{\vc^\prime(\mu^n)}  .  \eqno (5.2)$$
If $\r^n_s\tends\r_s$ in $L^1$ and there is $C>0$ such that 
$$\r_s^n(x)\le C\qquad\forall s\in[t,0]    ,
\txt{$m$ a. e.}x\in S   ,\quad\forall n\ge 1 \eqno (5.3)$$
then 
$$\ent_m(\mu_s)=
\lim_{n\tends+\infty}\ent_m(\mu^n_s)  \qquad
\forall s\in[t,0]   .  \eqno (5.4)$$

\proof The curves $\mu^n$ are absolutely continuous by proposition 2.4. Formula (5.1) follows from [3]; indeed, in [3] it is proven that, if 

\noindent a) the limit on the right hand side of (5.1) is finite (which in our case is true by point 3 of the hypotheses and proposition 2.4) and 

\noindent b) point 2) of the hypotheses hold, 

\noindent then $\mu$ is absolutely continuous and (5.1) holds. 

We prove (5.2). Since (2.10) defines $||\L_\mu||_{\vc^\prime(\mu)}$ as a $\sup$, it suffices to prove that, for all $\phi\in\tc$, the map
$$\fun{}{\mu}{
-\int_t^0\dr s\int_S[
\partial_s\phi_s+\frac{\b}{2}\D_\ec\phi_s
]  \dr\mu_s  +
\int_S\phi_0\dr\mu_0-\int_S\phi_t\dr\mu_t
}  $$
is continuous from $C([t,0],\pc(S))$ to $\R$. Since $\phi_0$ and $\phi_t$ are continuous functions, the last two terms on the right are continuous functions of $\mu$; thus, it suffices to show that
$$\fun{}{\mu}{
-\int_t^0\dr s\int_S[
\partial_s\phi_s+\frac{\b}{2}\D_\ec\phi_s
]  \dr\mu_s  
}  $$
is continuous. 

Let $\mu^n\tends\mu$ in 
$C([t,0],\pc(S))$; we have to show that 
$$\int_t^0\dr s\int_S[
\partial_s\phi_s+\frac{\b}{2}\D_\ec\phi_s
]\dr\mu^n_s\tends
\int_t^0\dr s\int_S[
\partial_s\phi_s+\frac{\b}{2}\D_\ec\phi_s
]\dr\mu_s  .  \eqno (5.5)$$
Since $\phi\in\tc$, we have that the integrand is in $C([t,0],L^\infty(S,m))$. Now it suffices to recall that the densities of $\mu^n_s$, being uniformly integrable, are weakly compact in 
$L^1$. 

We prove (5.4). Let $\{ \mu^{n^\prime}_s \}$ be a subsequence of $\{ \mu^n_s \}$; since 
$\r^n_s\tends\r_s$ in $L^1(m)$, we can find a subsequence $\{ \r^{n^\pprime}_s \}$ such that 
$\r^{n^\pprime}_s\tends\r_s$ $m$ a. e.. By (5.3) and dominated convergence, this implies that 
$$\int_S\r^{n^\pprime}_s\log\r^{n^\pprime}_s\dr m\tends
\int_S\r_s\log\r_s\dr m  .   $$
Thus, every subsequence of $\ent_m\mu^n_s$ has a sub-subsequence converging to 
$\ent_m\mu^n_s$, which implies (5.4).

\fin

We are going to prove formula (6) by an approximation argument; the first step is the next lemma, which shows that (6) holds for the "smoothened" curve $\mu^\e$ which we defined in section 3. To end the proof of theorem 1, we shall take advantage of the semicontinuity of lemma 5.1 and take limits as $\e$ goes to zero, getting (6). 

\lem{5.2} Let $t<0$ and let $\fun{\mu}{[t,0]}{\pc(S)}$ be an admissible curve such that (5) hold. Let us suppose, moreover, that $\mu$ is a weak solution of the Fokker-Planck equation with drift $V\in\vc(\mu)$. Let $\mu^\e=\r^\e m$ be defined as in section 3. Then, $\mu^\e$ satisfies the entropy generation equality (4.4).

\proof It suffices to check that $\mu^\e$ satisfies the hypotheses of lemma 4.2, i. e. that it is regular. First of all, point 1) of lemma 3.3 implies that $\r^\e$ is bounded. 

Next, point 2) of lemma 3.3 implies that there is $D_1(\e)>0$ such that, for all $s\in[t,0]$,
$$\int_S\G(\r^\e_s,\r^\e_s)\dr m\le D_1(\e)  .  $$
Together with point 1) of lemma 3.3 this implies that, for all $s\in[t,0]$, 
$$\int_S\frac{\G(\r^\e_s,\r^\e_s)}{\r^\e_s}\dr m\le D_2(\e)  .  \eqno (5.6)$$
Lastly, we have to show that $\fun{\mu^\e}{[t,0]}{\pc(S)}$ is Lipschitz; by lemma 3.3 we know that $\mu^\e$ solves the continuity equation for a drift $Y^\e\in\vc(\mu^\e)$. By proposition 2.4, Lipschitzianity follows if we show that there is $D_3(\e)$ such that, for all $s\in[t,0]$, 
$$||Y^\e_s||_{\zc(\mu^\e_s)}\le D_3(\e)   .   \eqno (5.7)$$
Let $\phi\in\tc$ and let $[a,b]\subset[t,0]$; by the formula before (2.16) we have the first equality below, while the second one is (1.5).  
$$\int_S\phi_b\dr\mu_b^\e-
\int_S\phi_a\dr\mu_a^\e-
\int_a^b\dr s\int_S\partial_s\phi_s\dr\mu^\e_s=$$
$$\frac{\b}{2}\int_a^b\dr s\int_S\D_\ec\phi_s\cdot\r^\e_s\dr m+
\int_a^b\inn{\phi_s}{V^\e_s}_{\zc(\mu^\e_s)}\dr s=$$
$$-\frac{\b}{2}\int_a^b\dr s\int_S\G(\phi_s,\r^\e_s)\dr m+
\int_a^b\inn{\phi_s}{V^\e_s}_{\zc(\mu^\e_s)}\dr s  .  $$
By Cauchy-Schwarz, this implies the inequality below.  
$$\left\vert
\int_S\phi_b\dr\mu_b^\e-
\int_S\phi_a\dr\mu_a^\e-
\int_a^b\dr s\int_S\partial_s\phi_s\dr\mu^\e_s
\right\vert   \le$$
$$\left[
\frac{\b}{2}\left(
\int_a^b\dr s\int_S\frac{\G(\r^\e_s,\r^\e_s)}{\r^\e_s}\dr m
\right)^\2    +
\left(
\int_a^b||V^\e_s||_{\zc(\mu^\e_s)}
\right)^\2
\right]  \cdot
\left(
\int_a^b\dr s\int_S\G(\phi_s,\phi_s)\r^\e_s\dr m
\right)^\2  .  \eqno (5.8)$$
Taking $[a,b]=[t,0]$, (5.8) implies by the Riesz representation theorem that $\mu^\e$ solves the continuity equation for a drift $Y^\e$ such that 
$$||Y^\e||_{\vc(\mu)}\le  
\left[
\frac{\b}{2}\left(
\int_t^0\dr s\int_S\frac{\G(\r^\e_s,\r^\e_s)}{\r^\e_s}\dr m
\right)^\2    +
\left(
\int_t^0||V^\e_s||_{\zc(\mu^\e_s)}
\right)^\2
\right]   .   $$
Note that the term on the right is finite by (5.6) and point 5) of lemma 3.3. 

Actually, since (5.8) holds for all intervals $[a,b]\subset[t,0]$, we easily see that, for a. e. 
$s\in[t,0]$, 
$$||Y^\e_s||_{\zc(\mu^\e_s)}\le
\left[
\frac{\b}{2}\left(
\int_S\frac{\G(\r^\e_s,\r^\e_s)}{\r^\e_s}\dr m
\right)^\2    +
\left(
||V^\e_s||_{\zc(\mu^\e_s)}
\right)^\2
\right]   .   $$
This implies (5.7) by (5.6) and point 5) of lemma 3.3.

\fin

\noindent{\bf Proof of theorem 1.} Let $\mu$ solve the Fokker-Planck equation with drift 
$V\in\vc(\mu)$; let (5) hold and let $\mu^\e=\r^\e m$ be defined as in section 3. By lemma 5.2 we have the equality below, where $O(\r^\e)$ is defined as in corollary 4.2; the inequality comes from point 4) of lemma 3.3. 
$$\inn{Y^\e}{Y^\e}_{\vc(\mu^\e)} +
\b\ent_m(\mu^\e_0)-\b\ent_m(\mu^\e_t)+
\frac{\b^2}{4}\inn{O(\r^\e)}{O(\r^\e)}_{\vc(\mu^\e)}=$$
$$||V^\e||^2_{\vc(\mu^\e)}\le e^{-K\b\e}||V||_{\vc(\mu)}  .  \eqno (5.9)$$
First of all, we show that $\mu\in AC^2([t,0],\pc(S))$. By point 1) of lemma 5.1 it suffices to show that 
$$\liminf_{\e\searrow 0}\inn{Y^\e}{Y^\e}_{\vc(\mu^\e)}<+\infty  .  $$
Actually, lemmas 5.1 and proposition 2.4 also imply that 
$$||Y||_{\vc(\mu)}\le
\liminf_{\e\searrow 0}\inn{Y^\e}{Y^\e}_{\vc(\mu^\e)}    \eqno (5.10)$$
provided we show that the right hand side is finite. In order to show this we note that, by (3.12), $\r^\e_{t-\e}=\hat\r^\e_{t-\e}$; in turn, by the way we extended $\hat\r_s$ outside 
$[t,0]$, we have $\hat\r_{t-\e}=\hat\r_t$; as $\e\tends 0$, 
$\hat\r^\e_{t}\tends\r_t$ in $L^1$ by (3.8). By lemma 5.1 this implies that 
$$\ent_m(\r^\e_{t-\e})\tends\ent_m(\r_t)  \txt{as}\e\tends 0  .  $$
Analogously, 
$$\ent_m(\r^\e_{\e})\tends\ent_m(\r_0)  \txt{as}\e\tends 0  .  $$
Now (6) follows from (5.10), the last two formulas and (5.9) on $[t-\e,\e]$.  

\fin

\vskip 2 pc
\centerline{\bf  \S 6}
\centerline{\bf Minimal characteristics}
\vskip 1 pc

\noindent{\bf Proof of corollary 1.} Point 1) comes from formula (5.2) of lemma 5.1. 

We prove point 2). Since (7) and (9) hold, each $\mu\in\fc$ satisfies the hypotheses of theorem 1; the inequality below comes from (6) and the definition of $M$ in (7). 
$$\int_t^0||\dot\mu_s||^2\dr s+
\b\ent_m(\mu_0)-\b\ent_m(\mu_t)\le M \qquad\forall\mu\in\fc .  $$
Since $\ent_m(\mu_0)\ge 0$ by Jensen's inequality, the last formula and (8) imply that there is $M^\prime>0$ such that  
$$\int_t^0||\dot\mu_s||^2\dr s\le M^\prime\qquad \forall\mu\in\fc  .  $$
By [3], this implies that $\fc$ is compact. 

We prove point 3). We begin to show that $U(t,\nu)$ is finite. First of all, the boundedness of $F$ and $U$ implies that 
$$U(t,\nu)>-\infty  .  $$
On the other side, let $\fun{\mu}{[t,0]}{\pc(S)}$ be the solution of the heat equation starting at $\nu$; clearly, $\mu$ solves the Fokker-Planck equation with drift $V\equiv 0$. Plugging this into (10) and using the fact that $F$ and $U$ are bounded, we get that
$$U(t,\nu)<+\infty   .   $$
Now we know that $U(t,\nu)\in\R$; let us prove that the $\inf$ in its definition is a minimum. 

We assert that there is a sequence $\mu^n=\r^nm$ and $\d_n\tends 0$ such that 
$$\d_n\le\r^n_s(x)\le\frac{1}{\d_n}
\txt{$m$-a. e.}x\in S\quad\forall s\in[t,0],   \eqno (6.1)$$
$$\r^n_t\tends \nu\txt{in $L^1(S,m)$ as $n\tends+\infty$}  \eqno (6.2)$$
and 
$$\liminf_{n\tends+\infty}
\2||\dot\mu^n-\frac{\b}{2}\D_\ec\mu^n||^2_{\vc^\prime(\mu^n)}+
\int_t^0\dr s\int_SF(s,x)\dr\mu^n_s(x)+
U(\mu^n_0)\le U(t,\mu_0)  .  \eqno (6.3)$$
To show this, we let $\tilde\mu^n$ be a sequence minimising in (10); we take 
$\e_n\searrow 0$ and we set $\mu^n_s=\hat\mu^{\e_n}_s$, where $\hat\mu^\e_s$ is defined as at the end of step 1 of section 3. Formula (6.1) follows by (1.9); (6.2) follows as in the proof of point 6) of lemma 3.3, (6.3) as in point 4) of the same lemma. 

Let us call $V^n$ the drift of the Fokker-Planck equation solved by $\mu^n$. Since 
$U(t,\nu)$ is finite and $\{ \mu^n \}$ is minimising, there is $M_4>0$ such that 
$$\2||V^n||^2_{\vc(\mu^n)}+
\int_t^0\dr s\int_SF(s,x)\dr\mu^n_s(x)+U(\mu^n_0)\le M_4  \eqno (6.3)$$
for all $n\ge 1$. Since $F$ and $U$ are bounded, this implies that
$$\2||V^n||^2_{\vc(\mu^n)}\le M_5\qquad\forall n\ge 1  .  \eqno (6.5)$$
By point 2), $\{ \mu^n \}$ converges, up to subsequences, to a curve $\mu$. 

We want to show that $\mu$ is a minimiser. This follows if we show that point 1) of this corollary holds, i. e. that $\mu^n_s=\r^n_sm$ with $\{ \r^n_s \}$ uniformly integrable. In the arguments before (2.13) we saw that 
$$||V^n||_{\vc(\mu^n)}=||\dot\mu^n-\frac{\b}{2}\D_\ec\mu_n||_{\vc^\prime(\mu^n)}  .  $$
Since $F$ and $U$ are bounded, (6.5) implies that 
$$||\dot\mu^n-\frac{\b}{2}\D_\ec\mu_n||_{\vc^\prime(\mu^n)}^2\le M_5\qquad
\forall n\ge 1  .  $$
By theorem 1, this implies that (6) holds, i. e. that, for all $s\in[t,0]$, 
$$\int_{t}^s||\dot\mu^n_s||^2\dr s+
\b\ent_m(\mu^n_s)-\b\ent_m(\mu^n_{t})\le M_5  .  \eqno (6.6)$$
By (6.2) and lemma 5.1 we note that $\ent(\mu^n_{t})$ is bounded as 
$\e_n\tends 0$; by (6.6) this implies that 
$$\ent_m(\mu^n_s)\le M_6\qquad
\forall s\in[t,0],\quad\forall n\ge 1   .  $$
Since the function $\fun{}{x}{x\log x}$ is superlinear, this implies that $\{ \r^n_s \}_{n,s}$ is equi-integrable. 

\fin

\vskip 2pc
\centerline{\bf Bibliography}





\noindent [1] L. Ambrosio, Lecture notes on optimal transport problems, in Mathematical aspects of evolving interfaces, LNS 1812, Springer, Berlin, 2003. 

\noindent [2] L. Ambrosio, N. Gigli, A. Mondino, T. Rajala, Riemannian Ricci curvature lower bounds in metric measure spaces with $\sigma$-finite measure, Trans. Amer. Math. Soc., {\bf 367}, 4661-4701, 2015.

\noindent [3] L. Ambrosio, N. Gigli, G. Savar\'e, Gradient Flows, Birkh\"auser, Basel, 2005.

\noindent [4] L. Ambrosio, N. Gigli, G. Savar\'e, Heat flow and calculus on metric measure spaces with Ricci curvature bounded below - the compact case. Analysis and numerics of Partial Differential Equations, 63-115, Springer, Milano, 2013.  

\noindent [5] L. Ambrosio, A. Mondino, G. Savar\'e, Nonlinear diffusion equations and curvature conditions in metric measure spaces, preprint 2017. 


\noindent [6] L. Ambrosio, N. Gigli, G. Savar\'e, Metric measure spaces with Riemannian Ricci curvature bounded from below, Duke Math. J., {\bf 163-7}, 1405-1490, 2014. 

\noindent [7] L. Ambrosio, N. Gigli, G. Savar\'e, Bakry-\'Emery curvature-dimension condition and Riemannian Ricci curvature bounds, Ann. Probab., {\bf 43}, 339-404, 2015.  







\noindent [8] U. Bessi, The stochastic value function in metric measure spaces, Discrete and Continuous Dynamical Systems, {\bf 37-4}, 1919-1839, 2017. 


\noindent [9] H. Brezis, Analisi Funzionale, Liguori, Napoli, 1986.


\noindent [10] J. Cheeger, Differentiability of Lipschitz functions on metric measure spaces, GAFA, Geom. Funct. Anal. {\bf 9}, 428-517, 1999.


\noindent [11] J. Feng, T. Nguyen, Hamilton-Jacobi equations in space of measures associated with a system of conservation laws, Journal de Math\'ematiques pures et Appliqu\'ees, {\bf 97}, 318-390, 2012. 




\noindent [12] W. H. Fleming, The Cauchy problem for a Nonlinear first order Partial Differential Equation, JDE, {\bf 5}, 515-530, 1969. 





\noindent [13] N. Gigli, B. Han, the continuity equation on metric measure spaces, Calc. Var. Partial Differential Equations, {\bf 53}, 149-177, 2015.

\noindent [14] S. Lisini, Characterisation of absolutely continuous curves in Wasserstein space, Calc. Var.Partial Differential Equations, {\bf 28}, 85-120, 2007.

\noindent [15] U. Mosco, Composite media and asymptotic Dirichlet forms, J. Functional Analysis, {\bf 123}, 368-421, 1994.





\noindent [16] E. Nelson, Dynamical theories of Brownian motion, Princeton, 1967.

\noindent [17] C. Villani, Topics in optimal transportation, Providence, R. I., 2003.

\end

By [7], the spaces $V^1_\infty$ and $V^2_\infty$ are $P_{-\e}$-invariant; thus, if 
$\phi\in\tc$, then $P_{-\e}(\phi)\in\tc$ too and for $\phi\in\tc$ we can define 
$$\tilde V^\e(\phi)=\inn{V}{P_{-\e}\phi}_{\vc(\mu)} ,\qquad
\tilde Y^\e(\phi)=\inn{Y}{P_{-\e}\phi}_{\vc(\mu)}   \eqno (3.10)$$
for the drifts $V$ of the Fokker-Planck and $Y$ of the continuity equation satisfied by 
$\mu$.

We are going to show that $\tilde V^\e$ and $\tilde Y^\e$ are bounded operators on $\tc$  for the norm $||\cdot||_{\vc(\tilde\mu^\e)}$; by Riesz's representation theorem, this will imply that they are represented by elements of $\vc(\tilde\mu^\e)$, which we continue to call, by an abuse of notation, $\tilde V^\e$ and $\tilde Y^\e$. 

In order to get a bound on the norm of $\tilde V^\e$ and $\tilde Y^\e$, we recall (1.14), which implies that
$$\G(P_{-\e}\phi,P_{-\e}\phi)\le e^{-K\b\e}P_{-\e}\G(\phi_s,\phi_s)$$
for the same $K$ in the definition of $RCD(K,\infty)$ space. This yields the second inequality below, while the first one is Cauchy-Schwarz in $\vc(\mu)$; the first equality below is (3.10); the second one comes from the duality between $P_{-\e}$ and $H_\e$ and the third one is the definition of $\tilde\mu^\e$ in (3.1). 
$$\tilde V^\e(\phi)=
\inn{V}{P_{-\e}\phi}_{\vc(\mu)}\le
||V||_{\vc(\mu)}\cdot
\left[
\int_t^0\dr s\int_S\G(P_{-\e}\phi_s,P_{-\e}\phi_s)\r_s\dr m
\right]^\2  \le$$
$$||V||_{\vc(\mu)}\cdot
\left[
\int_t^0\dr s\int_Se^{-K\b\e}P_{-\e}\G(\phi_s,\phi_s)\r_s\dr m
\right]^\2  =||V||_{\vc(\mu)}\cdot e^{-\frac{K\b}{2}\e}
\left[
\int_t^0\dr s\int_S\G(\phi_s,\phi_s)\dr (H_\e \mu_s)
\right]^\2  =  $$
$$||V||_{\vc(\mu)}\cdot e^{-\frac{K\b}{2}\e}
\left[
\int_t^0\dr s\int_S\G(\phi_s,\phi_s)\dr \tilde\mu_s^\e
\right]^\2   .   $$
From the last formula and (2.10) we get that 
$$||\tilde V^\e||_{\vc^\prime(\tilde\mu^\e)}\le e^{-\frac{K\b}{2}\e}
||V||_{\vc(\mu)}  .  \eqno (3.11)$$
Analogously, 
$$||\tilde Y^\e||_{\vc^\prime(\tilde\mu^\e)}\le e^{-\frac{K\b}{2}\e}
||Y||_{\vc(\mu)}  .  \eqno (3.12)$$